\documentclass[11pt]{amsart}

\usepackage{amssymb}
\usepackage{graphicx}
\usepackage[all]{xy}

\title{Comparing globular complex and flow}

\author{Philippe Gaucher}

\address{Preuves Programmes et Syst{\`e}mes\\ Universit{\'e} Paris 7--Denis Diderot\\
Case 7014\\2 Place Jussieu\\ 75251 PARIS Cedex 05\\ France}

\email{gaucher@pps.jussieu.fr}
\urladdr{http://www.pps.jussieu.fr/{\~{}}gaucher/}


\newcommand{\C}{\mathcal{C}}
\newcommand{\D}{\mathcal{D}}
\newcommand{\R}{\mathbb{R}}
\newcommand{\de}{\partial}
\newcommand{\p}\times
\renewcommand{\vec}{\overrightarrow}
\renewcommand{\P}{\mathbb{P}}

\newcommand{\beas}{\begin{eqnarray*}}
\newcommand{\eeas}{\end{eqnarray*}}

\newcounter{numerothm}[section]

\numberwithin{section}{part}

\newtheorem{thm}[numerothm]{Theorem}
\newtheorem{prop}[numerothm]{Proposition}
\newtheorem{lem}[numerothm]{Lemma}

\newtheorem{question}[numerothm]{Question}
\newtheorem{cor}[numerothm]{Corollary}

\newtheorem{defn}[numerothm]{Definition}
\newtheorem{propdef}[numerothm]{Proposition and Definition}

\newtheorem{nota}[numerothm]{Notation}
\newcommand{\bd}{\begin{defn}}
\newcommand{\ed}{\end{defn}}
\newcommand{\bp}{\begin{prop}}
\newcommand{\ep}{\end{prop}}
\newcommand{\bth}{\begin{thm}}
\renewcommand{\eth}{\end{thm}}
\newcommand{\bpf}{\begin{proof}}
\newcommand{\epf}{\end{proof}}

\newcommand{\fl}[1]{\ar@{->}[l]_{#1}}
\newcommand{\fr}[1]{\ar@{->}[r]^-{#1}}
\newcommand{\fd}[1]{\ar@{->}[d]_{#1}}
\newcommand{\fu}[1]{\ar@{->}[u]^{#1}}
\newcommand{\f}[2]{\ar@{->}[#1]|{#2}}
\newcommand{\ff}[2]{\ar@2{->}[#1]|{#2}}
\newcommand{\frr}[1]{\ar@{->}[rr]^{#1}}

\newcommand{\CW}{{\mathbf{CW}}}

\newcommand{\ho}{{\mathbf{Ho}}}
\newcommand{\iso}{\cong}

\newcommand{\vI}{\vec{I}}
\renewcommand{\leq}{\leqslant}
\renewcommand{\geq}{\geqslant}
\newcommand{\Rm}{\mathcal{R}^-}

\def\cartesien{%
  \ar@{-}[]+R+<6pt,-2pt>;[]+RD+<6pt,-6pt>%
  \ar@{-}[]+D+<2pt,-6pt>;[]+RD+<6pt,-6pt>%
}
\def\cocartesien{%
  \ar@{-}[]+L+<-6pt,+2pt>;[]+LU+<-6pt,+6pt>%
  \ar@{-}[]+U+<-2pt,+6pt>;[]+LU+<-6pt,+6pt>%
}

\newcommand{\brm}[1]{\rm{\mathbf{#1}}}

\renewcommand{\top}{{\brm{Top}}}
\newcommand{\gltop}{{\brm{glTop}}}

\newcommand{\glCW}{{\brm{glCW}}}
\newcommand{\dtop}{{\brm{Flow}}}

\newcommand{\set}{{\brm{Set}}}
\newcommand{\tdtop}{{\brm{FLOW}}}
\newcommand{\ttop}{{\brm{TOP}}}
\newcommand{\tgltop}{{\brm{glTOP}}}
\newcommand{\wdtop}{{\brm{qFlow}}}

\newcommand{\glob}{{\rm{Glob}}}

\newcommand{\liminj}{\varinjlim}
\newcommand{\limproj}{\varprojlim}

\DeclareMathOperator{\id}{Id}
\DeclareMathOperator{\Id}{Id}
\DeclareMathOperator{\card}{card}

\begin{document}

\begin{abstract}
A functor is constructed from the category of globular CW-comple\-xes
to that of flows. It allows the comparison of the S-homotopy
equivalences (resp. the T-homotopy equivalences) of globular complexes
with the S-homotopy equivalences (resp. the T-homotopy equivalences)
of flows.  Moreover, it is proved that this functor induces an
equivalence of categories from the localization of the category of
globular CW-complexes with respect to S-homotopy equivalences to the
localization of the category of flows with respect to weak S-homotopy
equivalences. As an application, we construct the underlying homotopy
type of a flow.
\end{abstract}

\maketitle

\tableofcontents

\part{Introduction}

\section{Outline of the paper}

The category of \textit{globular CW-complexes} $\glCW$ was introduced
in \cite{diCW} for modelling higher dimensional automata and
dihomotopy, the latter being an equivalence relation preserving their
computer-scientific properties, like the \textit{initial or final
states}, the presence or not of \textit{deadlocks} or of
\textit{unreachable states}, and more generally any computer-scientific 
property invariant by refinement of observation. More precisely, the
classes of \textit{S-homotopy equivalences} and of
\textit{T-homotopy equivalences} were defined. The category of
\textit{flows} as well as the notion of \textit{S-homotopy
equivalence} of flows are introduced in \cite{model3}. The notion of
S-homotopy equivalence of flows is interpreted in \cite{model3} as the
notion of homotopy arising from a model category structure.  The weak
equivalences of this model structure are called the
\textit{weak S-homotopy equivalences}.

The purpose of this paper is the comparison of the framework of
globular CW-complexes with the framework of flows. More precisely, we
are going to construct a functor $cat:\glCW\longrightarrow \dtop$ from
the category of globular CW-complexes to that of flows inducing an
equivalence of categories from the localization
$\glCW[\mathcal{SH}^{-1}]$ of the category of globular CW-complexes
with respect to the class $\mathcal{SH}$ of S-homotopy equivalences to
the localization $\dtop[\mathcal{S}^{-1}]$ of the category of flows
with respect to the class $\mathcal{S}$ of weak S-homotopy
equivalences.  Moreover, a class of T-homotopy equivalences of flows
will be constructed in this paper so that there exists, up to weak
S-homotopy, a T-homotopy equivalence of globular CW-complexes
$f:X\longrightarrow Y$ if and only if there exists a T-homotopy
equivalence of flows $g:cat(X)\longrightarrow cat(Y)$.

Part~\ref{part1} introduces the category of \textit{globular
complexes} $\gltop$, which is slightly larger than the category of
globular CW-complexes $\glCW$. Indeed, the latter category is not a
big enough setting for several constructions that are going to be
used.  Part~\ref{model1} builds the functor $cat:\gltop\longrightarrow
\dtop$. Part~\ref{model2} is a technical part which proves that two
globular complexes $X$ and $U$ are S-homotopy equivalent if and only
if the corresponding flows $cat(X)$ and $cat(U)$ are S-homotopy
equivalent.  Part~\ref{model4} proves that the functor
$cat:\glCW\longrightarrow \dtop$ from the category of globular
CW-complexes to that of flows induces an equivalence of categories
from the localization $\glCW[\mathcal{SH}^{-1}]$ of the category of
globular CW-complexes with respect to the class of S-homotopy
equivalences to the localization $\dtop[\mathcal{S}^{-1}]$ of the
category of flows with respect to the class of weak S-homotopy
equivalences. At last, Part~\ref{model2.5} studies and compares the
notion of T-homotopy equivalence for globular complexes and flows. And
Part~\ref{uhtf} applies all previous results to the construction of
the \textit{underlying homotopy type} of a flow.

\section{Warning}

This paper is the sequel of ``A model category for the homotopy theory
of concurrency'' \cite{model3}, where the category of flows was
introduced. This work is focused on the relation between the category
of globular CW-complexes and the category of flows. A first version of
the category of globular CW-complexes was introduced in a joined work
with Eric Goubault \cite{diCW}. A detailed abstract (in French) of
\cite{model3} and of this paper can be found in \cite{pgnote1} and
\cite{pgnote2}.

\section{Acknowledgment}

I thank the anonymous referee for the very careful reading of the
paper.

\part{S-homotopy and globular complex}\label{part1}

\section{Introduction}

The category of \textit{globular complexes} is introduced in
Section~\ref{def}. This requires the introduction of several other
notions, for instance the notion of \textit{multipointed topological
space}.  Section~\ref{stulim} carefully studies the behavior of the
functor $X\mapsto \tgltop(X,Y)$ for a given $Y$ with respect to the
globular decomposition of $X$ where $\tgltop(X,Y)$ is the set of
morphisms of globular complexes from $X$ to $Y$ equipped with the
Kelleyfication of the compact-open topology. At last,
Section~\ref{Shomgltop} defines and studies the notion of S-homotopy
equivalence of globular complexes. In particular, a cylinder functor
corresponding to this notion of equivalence is constructed.

\section{The category of globular complexes}\label{def}

\subsection{Compactly generated topological spaces}

The category $\top$ of \textit{compactly generated topological spaces}
(i.e. of weak Hausdorff $k$-spaces) is complete, cocomplete and
cartesian closed (more details for this kind of topological spaces in
\cite{MR90k:54001,MR2000h:55002}, the appendix of \cite{Ref_wH} and
also the preliminaries of \cite{model3}). Let us denote by
$\ttop(X,-)$ the right adjoint of the functor $-\p
X:\top\longrightarrow \top$. For any compactly generated topological
space $X$ and $Y$, the space $\ttop(X,Y)$ is the set of continuous
maps from $X$ to $Y$ equipped with the Kelleyfication of the
compact-open topology. For the sequel, any topological space will be
supposed to be compactly generated. A \textit{compact space} is always
Hausdorff.

\subsection{NDR pairs}

\bd Let $i:A\longrightarrow B$ and $p:X\longrightarrow Y$ be maps in a
category $\C$. Then $i$ has the \textit{left lifting property} (LLP)
with respect to $p$ (or $p$ has the \textit{right lifting property}
(RLP) with respect to $i$) if for any commutative square
\[
\xymatrix{
A\fd{i} \fr{\alpha} & X \fd{p} \\
B \ar@{-->}[ru]^{g}\fr{\beta} & Y}
\]
there exists $g$ making both triangles commutative. \ed

A Hurewicz fibration is a continuous map having the RLP with respect
to the continuous maps $\{0\}\p M\subset [0,1]\p M$ for any
topological space $M$. In particular, any continuous map having a
discrete codomain is a Hurewicz fibration. A Hurewicz cofibration is a
continuous map having the homotopy extension property. In the category
of compactly generated topological spaces, any Hurewicz cofibration is
a closed inclusion of topological spaces \cite{Ref_wH}. There exists a
model structure on the category of compactly generated topological
spaces such that the cofibrations are the Hurewicz cofibrations, the
fibrations are the Hurewicz fibrations, and the weak equivalences are
the homotopy equivalences (\cite{MR35:2284} \cite{MR39:4846}
\cite{ruse} and also \cite{strom2}). 
In this model structure, all topological spaces are fibrant and
cofibrant. The class of Hurewicz cofibrations coincides with the class
of NDR pairs. For any NDR pair $(Z,\de Z)$, one has
\cite{MR35:970}
\cite{MR80b:55001}
\cite{MR1802847} \cite{MR1867354}:
\begin{enumerate}
\item There exists a continuous map $\mu: Z\longrightarrow [0,1]$
such that $\mu^{-1}(\{0\})=\de Z$.
\item There exists
a continuous map $r: Z\p [0,1]\longrightarrow Z\p \{0\}
\cup \de Z\p [0,1]$ which is the identity on $ Z\p
\{0\} \cup \de Z\p [0,1]\subset  Z\p [0,1]$.
\end{enumerate}

This fact together with the continuous map $\mu: Z\longrightarrow
[0,1]$ is used in the proofs of Theorem~\ref{i} and of
Theorem~\ref{corT}.

\subsection{Definition of a globular complex}

A \textit{globular complex} is a topological space together with
a structure describing the sequential process of attaching
\textit{globular cells}. The class of globular complexes includes
the class of \textit{globular CW-complexes}. A general globular
complex may require an arbitrary long transfinite construction. We
must introduce this generalization because several constructions do
not stay within the class of globular CW-complexes.

\bd A {\rm multipointed topological space} $(X,X^0)$ is a pair of
topological spaces such that $X^0$ is a discrete subspace of $X$.  A
morphism of multipointed topological spaces $f:(X,X^0)\longrightarrow
(Y,Y^0)$ is a continuous map $f:X\longrightarrow Y$ such that
$f(X^0)\subset Y^0$. The corresponding category is denoted by
$\top^m$. The set $X^0$ is called the {\rm $0$-skeleton} of $(X,X^0)$.
The space $X$ is called the {\rm underlying topological space} of
$(X,X^0)$. \ed

A multipointed space of the form $(X^0,X^0)$ where $X^0$ is a
discrete topological space will be called a \textit{discrete
multipointed space} and will be frequently identified with $X^0$
itself.

\bp The category of multipointed topological spaces is
cocomplete. \ep

\bpf This is due to the facts that the category of topological
spaces is cocomplete and that the colimit of discrete spaces is a
discrete space. \epf

\bd Let $Z$ be a topological space. The {\rm globe of $Z$}, which
is denoted by $\glob^{top}(Z)$, is the multipointed space
\[(|\glob^{top}(Z)|,\{0,1\})\] where the topological space
$|\glob^{top}(Z)|$ is the quotient of $\{0,1\}\sqcup (Z\p[0,1])$
by the relations $(z,0)=(z',0)=0$ and  $(z,1)=(z',1)=1$ for any
$z,z'\in Z$. \ed

In particular, $\glob^{top}(\varnothing)$ is the multipointed
space $(\{0,1\},\{0,1\})$.

\begin{nota} If $Z$ is a singleton, then the
globe of $Z$ is denoted by $\vI^{top}$. \end{nota}

Any ordinal can be viewed as a small category whose objects are the
elements of $\lambda$, that is the ordinals $\gamma<\lambda$, and
where there exists a morphism $\gamma\longrightarrow \gamma'$ if and
only if $\gamma\leq \gamma'$.

\bd Let $\C$ be a cocomplete category. Let $\lambda$ be an
ordinal. A {\rm $\lambda$-sequence} in $\C$ is a colimit-preserving
functor $X:\lambda\longrightarrow \C$. Since $X$ preserves colimits,
for all limit ordinals $\gamma<\lambda$, the induced map
$\liminj_{\beta<\gamma}X_\beta\longrightarrow X_\gamma$ is an
isomorphism. The morphism $X_0\longrightarrow\liminj X$ is called the
{\rm transfinite composition} of $X$. \ed

\bd A {\rm relative globular precomplex} is a $\lambda$-sequence of
multipointed topological spaces $X:\lambda\longrightarrow \top^m$ such
that for any $\beta<\lambda$, there exists a pushout diagram of
multipointed topological spaces
\[\xymatrix{
\glob^{top}(\de Z_\beta)\fr{\phi_\beta}\fd{} & \fd{}X_\beta\\
\glob^{top}( Z_\beta)\fr{}& \cocartesien X_{\beta+1}}
\]
where the pair $(Z_\beta,\de Z_\beta)$ is a NDR pair of compact
spaces. The morphism \[\glob^{top}(\de Z_\beta)\longrightarrow
\glob^{top}( Z_\beta)\] is induced by the closed inclusion $\de
 Z_\beta\subset Z_\beta$. \ed

\bd A {\rm globular precomplex} is a $\lambda$-sequence of
multipointed topological spaces $X:\lambda\longrightarrow \top^m$ such
that $X$ is a relative globular precomplex and such that $X_0=(X^0,X^0)$
with $X^0$ a discrete space. \ed

Let $X$ be a globular precomplex. The $0$-skeleton of $\liminj X$ is
equal to $X^0$.

\bd A morphim of globular precomplexes $f:X\longrightarrow Y$ is a
morphism of multipointed spaces still denoted by $f$ from $\liminj X$
to $\liminj Y$. \ed

\begin{nota} If $X$ is a globular precomplex, then the underlying
topological space of the multipointed space $\liminj X$ is denoted by
$|X|$ and the $0$-skeleton of  the multipointed space $\liminj X$ is denoted by
$X^0$. \end{nota}

\bd Let $X$ be a globular precomplex. The space $|X|$ is called
the {\rm underlying topological space} of $X$. The set $X^0$ is called
the {\rm $0$-skeleton} of $X$. The family $(\de Z_\beta,
Z_\beta,\phi_\beta)_{\beta<\lambda}$ is called the {\rm globular
decomposition} of $X$. \ed

As set, the topological space $X$ is by construction the disjoint
union of $X^0$ and of the $|\glob^{top}( Z_\beta\backslash \de
 Z_\beta)|\backslash\{0,1\}$.

\bd Let $X$ be a globular precomplex. A morphism of globular
precomplexes $\gamma:\vI^{top}\longrightarrow X$ is a {\rm
non-constant execution path} of $X$ if there exists
$t_0=0<t_1<\dots<t_{n}=1$ such that:
\begin{enumerate}
\item  $\gamma(t_i)\in X^0$ for any $i$
\item $\gamma(]t_i,t_{i+1}[)\subset \glob^{top}(Z_{\beta_i}\backslash \de Z_{\beta_i})$ 
for some $(\de Z_{\beta_i},Z_{\beta_i})$ of the globular decomposition
of $X$
\item  for $0\leq i<n$, there exists $z^i_\gamma\in Z_{\beta_i}\backslash \de Z_{\beta_i}$
and a strictly increasing continuous map
$\psi^i_\gamma:[t_i,t_{i+1}]\longrightarrow [0,1]$ such that
$\psi^i_\gamma(t_i)=0$ and $\psi^i_\gamma(t_{i+1})=1$ and for any
$t\in [t_i,t_{i+1}]$, $\gamma(t)=(z^i_\gamma,\psi^i_\gamma(t))$.
\end{enumerate}
In particular, the  restriction
$\gamma\!\restriction_{]t_i,t_{i+1}[}$ of $\gamma$ to
$]t_i,t_{i+1}[$ is one-to-one. The set of non-constant execution
paths of $X$ is denoted by ${\P}^{ex}(X)$. \ed

\bd A morphism of globular precomplexes $f:X\longrightarrow Y$ is
{\rm non-decreasing} if the canonical set map
$\top([0,1],|X|)\longrightarrow \top([0,1],|Y|)$ induced by
composition by $f$ yields a set map ${\P}^{ex}(X)\longrightarrow
{\P}^{ex}(Y)$. In other terms, one has the commutative diagram of
sets
\[\xymatrix{
{\P}^{ex}(X)\fr{}\fd{\subset}& {\P}^{ex}(Y)\fd{\subset}\\
\top([0,1],|X|) \fr{} &\top([0,1],|Y|)}
\]
\ed

\bd A {\rm globular complex } (resp. a {\rm relative globular complex})
$X$ is a globular precomplex (resp. a relative globular precomplex)
such that the attaching maps $\phi_\beta$ are non-decreasing. A
morphism of globular complexes is a morphism of globular precomplexes
which is non-decreasing. The category of globular complexes together
with the morphisms of globular complexes as defined above is denoted
by $\gltop$. The set $\gltop(X,Y)$ of morphisms of globular complexes
from $X$ to $Y$ equipped with the Kelleyfication of the compact-open
topology is denoted by $\tgltop(X,Y)$. \ed

Forcing the restrictions $\gamma\!\restriction_{]t_i,t_{i+1}[}$ to be
one-to-one means that only the ``stretched situation'' is
considered. It would be possible to build a theory of non-stretched
execution paths, non-stretched globular complexes and non-stretched
morphisms of globular complexes but this would be without interest
regarding the complexity of the technical difficulties we would meet.

\bd Let $X$ be a globular complex. A point $\alpha$ of $X^0$ such that
there are no non-constant execution paths ending to $\alpha$
(resp. starting from $\alpha$) is called {\rm initial state} (resp.
{\rm final state}). More generally, a point of $X^0$ will be sometime
called {\rm a state} as well. \ed

A very simple example of globular complex is obtained by
concatenating globular complexes of the form
$\glob^{top}(Z_j)$ for $1\leq i \leq n$ by identifying
the final state $1$ of
$\glob^{top}(Z_j)$ with the initial state
$0$ of $\glob^{top}(Z_{j+1})$.

\begin{nota} This
globular complex will be denoted by
\[\glob^{top}(Z_1)*\glob^{top}(Z_2)*\dots * \glob^{top}(Z_n)\]
\end{nota}

\begin{figure}
\begin{center}
\includegraphics[width=7cm]{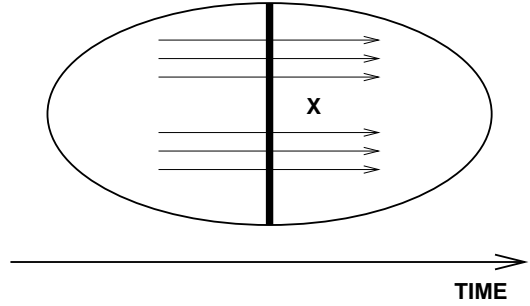}
\end{center}
\caption{Symbolic representation of
$\glob^{top}(X)$ for some compact topological space
$X$} \label{exglob}
\end{figure}

\subsection{Globular CW-complex}

Let $n\geq 1$. Let $\mathbf{D}^n$ be the closed $n$-dimensional disk
defined by the set of points $(x_1,\dots,x_n)$ of $\R^n$ such that
$x_1^2+\dots +x_n^2\leq 1$ endowed with the topology induced by that
of $\R^n$. Let $\mathbf{S}^{n-1}=\de \mathbf{D}^n$ be the boundary of
$\mathbf{D}^n$ for $n\geq 1$, that is to say the set of
$(x_1,\dots,x_n)\in \mathbf{D}^n$ such that $x_1^2+\dots
+x_n^2=1$. Notice that $\mathbf{S}^0$ is the discrete two-point
topological space $\{-1,+1\}$.  Let $\mathbf{D}^0$ be the one-point
topological space. Let $\mathbf{S}^{-1}$ be the empty space.

\bd \cite{diCW}\label{defCW} A {\rm globular
CW-complex} $X$ is a globular complex such that its globular
decomposition $(\de Z_\beta, Z_\beta,\phi_\beta)_{\beta<\lambda}$
satisfies the following properties. There exists a strictly increasing
sequence $(\kappa_n)_{n\geq 0}$ of ordinals with $\kappa_0=0$,
$\sup_{n\geq 0}\kappa_n=\lambda$, and such that for any $n\geq 0$, one
has the following fact:
\begin{enumerate}
\item for any $\beta\in[\kappa_n,\kappa_{n+1}[$, 
$(Z_\beta,\de Z_\beta)=(\mathbf{D}^{n},\mathbf{S}^{n-1})$
\item one has the pushout of multipointed topological spaces
\[
\xymatrix{ 
\bigsqcup_{i\in [\kappa_n,\kappa_{n+1}[}
\glob^{top}(\mathbf{S}^{n-1})\fr{\phi_n}\fd{} & X_{\kappa_n}\fd{}\\
\bigsqcup_{i\in [\kappa_n,\kappa_{n+1}[} \glob^{top}(\mathbf{D}^{n}) \fr{} &
X_{\kappa_{n+1}}\cocartesien}
\]
\end{enumerate}
where $\phi_n$ is the morphism of globular complexes induced by the
$\phi_\beta$ for $\beta\in[\kappa_n,\kappa_{n+1}[$. The full and
faithful subcategory of $\gltop$ of globular CW-complexes is denoted
by $\glCW$. Notice that we necessarily have $\liminj_n X_{\kappa_n}=X$ 
\ed

One also has:

\bp\cite{diCW}\label{globeCW} The Globe functor $X\mapsto
\glob^{top}(X)$ induces a functor from the category of
CW-complexes to the category of globular
CW-complexes. \ep

\section{Morphisms of globular complexes and colimits}\label{stulim}

The category of general topological spaces is denoted by
$\mathcal{T}$.

\bp\label{incltoptop} The inclusion of sets
$i:\tgltop(X,Y)\longrightarrow \ttop(|X|,|Y|)$ is an
inclusion of topological spaces, that is
$\tgltop(X,Y)$ is the subset of morphisms of globular
complexes of the space $\ttop(|X|,|Y|)$ equipped with
the Kelleyfication of the relative topology. \ep

\bpf Let $\rm{Cop}(|X|,|Y|)$ be the set of continuous maps
from $|X|$ to $|Y|$ equipped with the compact-open topology. The
continuous map \[\gltop(X,Y)\cap
\rm{Cop}(|X|,|Y|)\longrightarrow \rm{Cop}(|X|,|Y|)\] is an
inclusion of topological spaces. Let $f:Z\to
k(\rm{Cop}(|X|,|Y|))$ be a continuous map such
that $f(Z)\subset \gltop(X,Y)$ where $Z$ is
an object of $\top$ and where $k(-)$ is the Kelleyfication
functor. Then $f:Z\longrightarrow \rm{Cop}(|X|,|Y|)$ is continuous
since the Kelleyfication is a right adjoint and since $Z$ is a
$k$-space. So $f$ induces a continuous map $Z\to
\gltop(X,Y)\cap \rm{Cop}(|X|,|Y|)$, and therefore a
continuous map \[Z\longrightarrow k(\gltop(X,Y)\cap
\rm{Cop}(|X|,|Y|))\iso \tgltop(X,Y).\] \epf

\bp \label{limrelativetop} Let $(X_i)$ and $(Y_i)$ be two diagrams of
objects of $\mathcal{T}$. Let $f:(X_i)\longrightarrow (Y_i)$ be a
morphism of diagrams such that for any $i$, $f_i:X_i\longrightarrow
Y_i$ is an inclusion of topological spaces, i.e. $f_i$ is one-to-one
and $X_i$ is homeomorphic to $f(X_i)$ equipped with the relative
topology coming from the set inclusion $f(X_i)\subset Y_i$. Then the
continuous map $\limproj X_i\longrightarrow \limproj Y_i$ is an
inclusion of topological spaces, the limits $\limproj X_i$ and
$\limproj Y_i$ being calculated in $\mathcal{T}$.  \ep

Loosely speaking, the lemma above means that the limit in
$\mathcal{T}$ of the relative topology is the relative topology
of the limit.

\bpf Saying that $X_i\longrightarrow Y_i$ is an inclusion of
topological spaces is equivalent to saying that the isomorphism
of sets
\[\mathcal{T}(Z,X_i)\iso \{f\in \mathcal{T}(Z,Y_i);f(X_i)\subset
Y_i\}\] holds for any $i$ and for any object $Z$ of $\mathcal{T}$.
But like in any category, one has the isomorphism of sets
\[\limproj \mathcal{T}(Z,X_i)\iso \mathcal{T}(Z,\limproj X_i)\]
and
\[\limproj \mathcal{T}(Z,Y_i)\iso \mathcal{T}(Z,\limproj Y_i).\]
Using the construction of limits in the category of sets, it is
then obvious that the set $\mathcal{T}(Z,\limproj X_i)$ is
isomorphic to the set
\[\{f\in \limproj \mathcal{T}(Z,Y_i);f_i(X_i)\in Y_i\}\]
for any object $Z$ of $\mathcal{T}$. Hence the result. \epf

\bth\label{limittop} Let $X$ be a globular complex with globular
decomposition \[(\de Z_\beta, Z_\beta,\phi_\beta)_{\beta<\lambda}.\]
Then for any limit ordinal $\beta\leq \lambda$, one has the
homeomorphism \[\tgltop(X_\beta,U)\iso \limproj_{\alpha<\beta}
\tgltop(X_\alpha,U).\] And for any $\beta<\lambda$, one has the
pullback of topological spaces
\[
\xymatrix{\tgltop(X_{\beta+1},U)\fr{} \fd{}\cartesien &
\tgltop(\glob^{top}( Z_\beta),U)\fd{}\\
\tgltop(X_{\beta},U)\fr{} & \tgltop(\glob^{top}(\de Z_\beta),U)}
\]
\eth

\bpf One has the isomorphism of sets
\[\gltop(X_\beta,U)\iso \limproj_{\alpha<\beta}
\gltop(X_\alpha,U)\] and the pullback of sets
\[
\xymatrix{\gltop(X_{\beta+1},U)\fr{} \fd{}\cartesien &
\gltop(\glob^{top}( Z_\beta),U)\fd{}\\
\gltop(X_{\beta},U)\fr{} & \gltop(\glob^{top}(\de Z_\beta),U)}
\]
One also has the isomorphism of topological spaces
\[\ttop(|X_\beta|,|U|)\iso \limproj_{\alpha<\beta}
\ttop(|X_\alpha|,|U|)\] and the pullback of spaces
\[
\xymatrix{\ttop(|X_{\beta+1}|,|U|)\fr{} \fd{}\cartesien &
\ttop(|\glob^{top}( Z_\beta)|,|U|)\fd{}\\
\ttop(|X_{\beta}|,|U|)\fr{} & \ttop(|\glob^{top}(\de
 Z_\beta)|,|U|)}
\]
The theorem is then a consequence of
Proposition~\ref{limrelativetop}, of Proposition~\ref{incltoptop}
and of the fact that the Kelleyfication functor is a right
adjoint which therefore preserves all limits. \epf

\bp\label{decompdiscret} Let $X$ and $U$ be two globular
complexes. Then one has the homeomorphism 
\[\tgltop(X,U)\iso\bigsqcup_{\phi:X^0\longrightarrow U^0} \{f\in
\tgltop(X,U),f^0=\phi\}.\] 
\ep

\bpf The composite set map \[\tgltop(X,U)\to
\ttop(X,U)\longrightarrow \ttop(X^0,U^0)\] is
continuous and $\ttop(X^0,U^0)$ is a discrete
topological space. \epf

Let $X$ be a globular complex. The set ${\P}^{ex}X$ of non-constant
execution paths of $X$ can be equipped with the Kelleyfication of the
compact-open topology.  The mapping ${\P}^{ex}$ yields a functor from
$\gltop$ to $\top$ by sending a morphism of globular complexes $f$ to
$\gamma \mapsto f\circ
\gamma$.

\bd A {\rm globular subcomplex} $X$ of a globular complex $Y$ is
a globular complex $X$ such that the underlying topological space is
included in the one of $Y$ and such that the inclusion map $X\subset
Y$ is a morphism of globular complexes. \ed

\bp\label{pathtopS} Let $X$ be a globular complex. 
Then there is a natural isomorphism of topological spaces
$\tgltop(\vec{I}^{top},X)\iso {\P}^{ex} X$. \ep

\bpf Obvious. \epf

\bp\label{decomp} Let $Z$ be a topological space. Then one has
the isomorphism of topological spaces ${\P}^{ex}(\glob^{top}(Z))\iso
Z\p\tgltop(\vI^{top},\vI^{top})$. \ep

\bpf There is a canonical inclusion
\[{\P}^{ex}(\glob^{top}(Z))\subset
\ttop([0,1],Z\p [0,1]).\] The image of this inclusion
is exactly the subspace of \[f=(f_1,f_2)\in
\ttop([0,1],Z\p [0,1])\] such that
$f_1:[0,1]\longrightarrow Z$ is a constant map and such that
$f_2:[0,1]\longrightarrow [0,1]$ is a non-decreasing continuous map
with $f_2(0)=0$ and $f_2(1)=1$. Hence the
isomorphism of topological spaces. \epf

\section{S-homotopy in $\gltop$}\label{Shomgltop}

\subsection{S-homotopy in $\gltop$}

We now recall the notion of S-homotopy introduced in \cite{diCW} for a
particular case of globular complex.

\bd Two morphisms of globular complexes $f$ and $g$ from $X$ to
$Y$ are said {\rm S-homotopic} or {\rm S-homotopy equivalent} if there
exists a continuous map $H:[0,1]\p X\longrightarrow Y$ such that for
any $u\in [0,1]$, $H_u=H(u,-)$ is a morphism of globular complexes
from $X$ to $Y$ and such that $H_0=f$ and $H_1=g$. We denote this
situation by $f\sim_{S} g$. \ed

Proposition~\ref{pathtopS} justifies the following definition.

\bd Two execution paths of a globular complex $X$ 
are {\rm S-homoto\-pic} or {\rm S-homotopy equivalent} if the
corresponding morphisms of globular complexes from $\vI^{top}$ to $X$
are S-homoto\-py equivalent. \ed

\bd Two globular complexes $X$ and $Y$ are {\rm S-homotopy equivalent}
if and only if there exists two morphisms of $\gltop$
$f:X\longrightarrow Y$ and $g:Y\longrightarrow X$ such that $f\circ
g\sim_{S} Id_Y$ and $g\circ f\sim_{S} Id_X$.  This defines an
equivalence relation on the set of morphisms between two given
globular complexes called S-homotopy. The maps $f$ and $g$ are called
{\rm S-homotopy equivalence}. The mapping $g$ is called a {\rm
S-homotopic inverse of $f$}. \ed

\subsection{Pairing $\boxtimes$ between a compact topological space and a globular complex}

Let $U$ be a compact topological space.  Let $X$ be a globular complex
with the globular decomposition $(\de Z_\beta,
Z_\beta,\phi_\beta)_{\beta<\lambda}$. Let $(U\boxtimes
X)_0:=(X^0,X^0)$. If $Z$ is any topological space, let $U\boxtimes
\glob^{top}(Z):=\glob^{top}(U\p Z)$.

If $(Z,\de Z)$ is a NDR pair, then the continuous map $i:[0,1]\p \de
Z\cup \{0\}\p Z\longrightarrow [0,1]\p Z$ has a retract $r:[0,1]\p
Z\longrightarrow [0,1]\p \de Z\cup \{0\}\p Z$. Therefore $i\p
\id_U:[0,1]\p \de Z\p U\cup \{0\}\p Z\p U\longrightarrow [0,1]\p Z\p
U$ has a retract $r\p \id_U:[0,1]\p Z\p U\longrightarrow [0,1]\p \de
Z\p U\cup \{0\}\p Z\p U$. Therefore $(U\p Z,U\p\de Z)$ is a NDR pair.

Let us suppose $(U\boxtimes X)_\beta$ defined for an ordinal $\beta$
such that $\beta+1<\lambda$ and assume that $(U\boxtimes X)_\beta$ has
the globular decomposition $(U\p\de Z_\mu,U\p
Z_\mu,\psi_\mu)_{\mu<\beta}$. From the morphism of globular complexes
$\phi_\beta:\glob^{top}(\de Z_\beta)\longrightarrow X_\beta$, one
obtains the morphism of globular complexes
$\psi_\beta:\glob^{top}(U\p\de Z_\beta)\longrightarrow (U\boxtimes
X)_\beta$ defined as follows: an element $\phi_\beta(z)$ belongs to a
unique $Z_\mu \backslash \de Z_\mu$.  Then let
$\psi_\beta(u,z)=(u,\phi_\beta(z))$. Then let us define $(U\boxtimes
X)_{\beta+1}$ by the pushout of multipointed topological spaces
\[\xymatrix{
U\boxtimes\glob^{top}(\de Z_\beta)\fd{}\fr{\psi_\beta}& U\boxtimes X_\beta \fd{}\\
U\boxtimes\glob^{top}( Z_\beta)\fr{} & \cocartesien U\boxtimes X_{\beta+1}}
\]
Then the globular decomposition of  $(U\boxtimes X)_{\beta+1}$ is
$(U\p\de Z_\mu,U\p Z_\mu,\psi_\mu)_{\mu<\beta+1}$. If $\beta\leq
\lambda$ is a limit ordinal, let $(U\boxtimes
X)_\beta=\liminj_{\mu<\beta}(U\boxtimes X)_\mu$ as multipointed
topological spaces.

\bp\label{boxq} Let $U$ be a compact space. Let $X$ be a globular
complex. Then the underlying space $|U\boxtimes X|$ of $U\boxtimes X$
is homeomorphic to the quotient of $U\p |X|$ by the equivalence
relation making the identification $(u,x)=(u',x)$ for any $u,u'\in U$
and for any $x\in X^0$ and equipped with the final topology. \ep

\bpf The graph of this equivalence relation is $\Delta U\p |X| \p |X| \subset
U\p U\p |X| \p |X|$ where $\Delta U$ is the diagonal of $U$. It is a
closed subspace of $U\p U\p |X| \p |X|$. Therefore the quotient set
equipped with the final topology is still weak Hausdorff, and
therefore compactly generated. It then suffices to proceed by
transfinite induction on the globular decomposition of $X$. \epf

The underlying set of $U\boxtimes X$ is then exactly equal to
$X^0\sqcup (U\p (X\backslash X^0))$.  The point $(u,x)$ with $x\in
X\backslash X^0$ will be denoted also by $u\boxtimes x$.  If $x\in
X^0$, then by convention $u\boxtimes x=u'\boxtimes x$ for any $u,u'\in
[0,1]$.

\bp Let $U$ and $V$ be two compact spaces. Let $X$ be a globular complex.
Then there exists a natural morphism of globular complexes $(U\p
V)\boxtimes X\iso U\boxtimes (V\boxtimes X)$. \ep

\bpf Transfinite induction on the globular decomposition of $X$. \epf

\subsection{Cylinder functor for S-homotopy in $\gltop$}

\bp \label{caracgl} Let $f$ and $g$ be two morphisms of globular
complexes from $X$ to $Y$. Then $f$ and $g$ are S-homotopic if
and only if there exists a continuous map
\[h\in\top([0,1],\tgltop(X,Y))\] such that
$h(0)=f$ and $h(1)=g$. \ep

\bpf Suppose that $f$ and $g$ are S-homotopic. Then the S-homotopy $H$
yields a continuous map \[h\in\top([0,1]\p |X|,|Y|)\iso
\top([0,1],\ttop(|X|,|Y|))\] by construction, and $h$ is necessarily
in \[\top([0,1],\linebreak[0]\tgltop(X,Y))\] by hypothesis.
Conversely, if \[h\in \top([0,1],\tgltop(X,Y))\] is such that $h(0)=f$
and $h(1)=g$, then the isomorphism \[ \top([0,1],\ttop(|X|,|Y|))\iso
\top([0,1]\p |X|,|Y|)\] provides a map $H\in \top([0,1]\p |X|,|Y|)$
which is a S-homotopy from $f$ to $g$. \epf

\bth\label{existencegl} Let $U$ be a connected non-empty topological
space.  Let $X$ and $Y$ be two globular complexes. Then there exists
an isomorphism of sets \[\gltop(U\boxtimes X,Y)\iso
\top(U,\tgltop(X,Y)).\] \eth

\bpf If $X$ is a singleton (this implies in particular that $X=X^0$),
then $U\boxtimes X=X$. So in this case, $\gltop(U\boxtimes X,Y)\iso
\top(U,\tgltop(X,Y))\iso Y^0$ \textit{since $U$ is connected and
non-empty} and by Proposition~\ref{decompdiscret}. Now if
$X=\glob^{top}(Z)$ for some compact space $Z$, then
\[\gltop(U\boxtimes X,Y)\iso \gltop(\glob^{top}(Z\p U),Y)\] and it is
straightforward to check that the latter space is isomorphic to
\[\top(U,\tgltop(\glob^{top}(Z),Y)).\] Hence the isomorphism
\[\gltop(U\boxtimes X,Y)\iso \top(U,\tgltop(X,Y))\] if $X$ is a point
or a globe.

Let $(\de Z_\beta, Z_\beta,\phi_\beta)_{\beta<\lambda}$ be the
globular decomposition of $X$. Then one deduces that
\[\gltop((U\boxtimes X)_\beta,Y)\iso \top(U,\tgltop(X_\beta,Y)).\] for
any $\beta$ by an easy transfinite induction, using the construction
of $U\boxtimes X$ and Theorem~\ref{limittop}.  \epf

\bd Let $\C$ be a category. A {\rm cylinder} is a functor
$I:\C\longrightarrow \C$ together with natural transformations
$i_0,i_1:\id_\C\longrightarrow I$ and $p:I\longrightarrow \id_\C$ such that
$p\circ i_0$ and $p\circ i_1$ are the identity natural
transformation. \ed

\begin{cor}\label{cylgl}
The mapping $X\mapsto [0,1]\boxtimes X$ induces a functor from
$\gltop$ to itself which is a cylinder functor with the natural
transformations $e_i:\{i\}\boxtimes -
\longrightarrow [0,1]\boxtimes -$ induced by the inclusion maps
$\{i\}\subset [0,1]$ for $i\in\{0,1\}$ and with the natural
transformation $p:[0,1]\boxtimes -\longrightarrow \{0\}\boxtimes -$
induced by the constant map $[0,1]\longrightarrow \{0\}$. Moreover,
two morphisms of globular complexes $f$ and $g$ from $X$ to $Y$
are S-homotopic if and only if there exists a morphism of
globular complexes $H:[0,1]\boxtimes X\longrightarrow Y$ such that
$H\circ e_0=f$ and $H\circ e_1=g$. Moreover $e_0\circ H\sim_{S}
Id$ and $e_1\circ H\sim_{S} Id$.
\end{cor}

\bpf Consequence of Proposition~\ref{caracgl} and
Theorem~\ref{existencegl}. \epf

\section{Conclusion}

We are now ready for the construction of the functor
$cat:\gltop\longrightarrow \dtop$.

\part{Associating a flow with any globular CW-complex}\label{model1}

\section{Introduction}

After a short reminder about the category of flows in
Section~\ref{reminderflow}, the functor $cat:\gltop\longrightarrow
\dtop$ is constructed in Section~\ref{embedding}. For that
purpose, the notion of \textit{quasi-flow} is introduced.
Section~\ref{expush} comes back to the case of flows by
explicitely calculating the pushout of a morphism of flows of the
form $\glob(\de Z)\longrightarrow \glob(Z)$. This will be used in
Section~\ref{stretched} and in Part~\ref{model4}.
Section~\ref{stretched} proves that for any globular complex $X$,
the natural continuous map ${\P}^{top}X\longrightarrow {\P} X$
has a right hand inverse $i_X:{\P} X \longrightarrow {\P}^{top}X$
(Theorem~\ref{i}). The latter map has no reason to be natural.

\section{The category of flows}\label{reminderflow}

\bd \cite{model3} A {\rm flow} $X$ consists of a topological space
${\P} X$, a discrete space $X^0$, two continuous maps $s$ and $t$
from ${\P} X$ to $X^0$ and a continuous and associative map
$*:\{(x,y)\in {\P} X\p {\P} X; t(x)=s(y)\}\longrightarrow {\P} X$
such that $s(x*y)=s(x)$ and $t(x*y)=t(y)$.  A morphism of flows
$f:X\longrightarrow Y$ consists of a set map
$f^0:X^0\longrightarrow Y^0$ together with a continuous map ${\P}
f:{\P} X\longrightarrow {\P} Y$ such that $f(s(x))=s(f(x))$,
$f(t(x))=t(f(x))$ and $f(x*y)=f(x)*f(y)$. The corresponding
category is denoted by $\dtop$. \ed

The continuous map $s:{\P} X\longrightarrow X^0$ is called the
\textit{source map}. The continuous map $t:{\P} X\longrightarrow X^0$
is called the \textit{target map}. One can canonically extend
these two maps to the whole underlying topological space
$X^0\sqcup {\P} X$ of $X$ by setting
$s(x)=x$ and $t(x)=x$ for $x\in X^0$.

The topological space $X^0$ is called the \textit{$0$-skeleton of
$X$}\footnote{The reason of this terminology: the $0$-skeleton of a
flow will correspond to the $0$-skeleton of a globular CW-complex by
the functor $cat$ ; one could define for any $n\geq 1$ the
$n$-skeleton of a globular CW-complex in an obvious way.}. The
$0$-dimensional elements of $X$ are called \textit{states} or
\textit{constant execution path}.

The elements of ${\P} X$ are called \textit{non-constant execution
path}.  If $\gamma_1$ and $\gamma_2$ are two non-constant execution
paths, then $\gamma_1 *\gamma_2$ is called the concatenation or the
composition of $\gamma_1$ and $\gamma_2$. For $\gamma\in {\P} X$,
$s(\gamma)$ is called the \textit{beginning} of $\gamma$ and
$t(\gamma)$ the \textit{ending} of $\gamma$.

\begin{nota} 
For $\alpha,\beta\in X^0$, let ${\P}_{\alpha,\beta}X$ be the subspace
of ${\P} X$ equipped the Kelleyfication of the relative topology
consisting of the non-constant execution paths of $X$ with beginning
$\alpha$ and with ending $\beta$.
\end{nota}

\bd\cite{model3} 
Let $Z$ be a topological space. Then the {\rm globe} of $Z$ is the
flow $\glob(Z)$ defined as follows: $\glob(Z)^0=\{0,1\}$,
${\P}\glob(Z)=Z$, $s(z)=0$, $t(z)=1$ for any $z\in Z$ and the
composition law is trivial. 
\ed

\bd\cite{model3} 
The {\rm directed segment} $\vI$ is the flow defined as follows:
$\vI^0=\{0,1\}$, ${\P}\vI=\{[0,1]\}$, $s=0$ and $t=1$.
\ed

\bd Let $X$ be a flow. A point $\alpha$ of $X^0$ such that there are
no non-constant execution paths $\gamma$ such that $t(\gamma)=\alpha$
(resp.  $s(\gamma)=\alpha$) is called \textit{initial state}
(resp. \textit{final state}).  \ed

\begin{nota} 
The space $\tdtop(X,Y)$ is the set $\dtop(X,Y)$ equipped with the
Kelleyfication of the compact-open topology.
\end{nota}

\bp (\cite{model3} Proposition~4.15) \label{can} 
Let $X$ be a flow. Then one has the following natural isomorphism of
topological spaces ${\P} X\iso\tdtop(\vI,X)$.
\ep

\bth (\cite{model3} Theorem~4.17) 
\label{lim-colim} The category $\dtop$ is complete
and cocomplete. In particular, a terminal object is the flow
$\mathbf{1}$ having the discrete set $\{0,u\}$ as underlying
topological space with $0$-skeleton $\{0\}$ and with path space
$\{u\}$.  And the initial object is the unique flow $\varnothing$ having
the empty set as underlying topological space.  \eth

\bth \label{commute} (\cite{model3} Theorem~5.10) 
The mapping
\[(X,Y)\mapsto \tdtop(X,Y)\] induces a functor from $\dtop\p \dtop$ to
$\top$ which is contravariant with respect to $X$ and covariant with
respect to $Y$. Moreover: \begin{enumerate} \item One has the
homeomorphism \[\tdtop(\liminj_i X_i,Y)\iso \limproj_i \tdtop(X_i,Y)\]
for any colimit $\liminj_i X_i$ in $\dtop$.  \item For any finite
limit $\limproj_i X_i$ in $\dtop$, one has the homeomorphism
\[\tdtop(X,\limproj_i Y_i)\iso \limproj_i\tdtop(X, Y_i).\]
\end{enumerate} \eth

\section{The functor $cat$ from  $\gltop$ to $\dtop$} \label{embedding}

The purpose of this section is the proof of the following theorems:

\bth There exists a unique functor
$cat:\gltop\longrightarrow\dtop$ such that
\begin{enumerate}
\item if $X=X^0$ is a discrete globular complex, then $cat(X)$ is
the achronal flow $X^0$ (``achronal'' meaning with an empty path space)
\item for any compact topological space $Z$, $cat(\glob^{top}(Z))=\glob(Z)$
\item for any globular complex $X$ with globular decomposition
$(\de Z_\beta, Z_\beta,\phi_\beta)_{\beta<\lambda}$, for any limit
ordinal $\beta\leq\lambda$, the canonical morphism of flows
\[\liminj_{\alpha<\beta} cat(X_\alpha)\longrightarrow
cat(X_\beta)\] is an isomorphism of flows
\item for any globular complex $X$ with globular decomposition
$(\de Z_\beta, Z_\beta,\phi_\beta)_{\beta<\lambda}$, for any
$\beta<\lambda$, one has the pushout of flows
\[\xymatrix{\glob(\de Z_\beta)\fr{cat(\phi_\beta)}\fd{}& cat(X_\beta)\fd{}\\
\glob( Z_\beta)\fr{} & cat(X_{\beta+1})\cocartesien}\]
\end{enumerate}
\eth

\begin{nota} Let $M$ be a topological space. Let $\gamma_1$ and
$\gamma_2$ be two continuous maps from $[0,1]$ to $M$ with
$\gamma_1(1)=\gamma_2(0)$. Let us denote by
$\gamma_1 *_a \gamma_2$ (with $0<a<1$) the following continuous
map: if $0\leq t\leq a$, $(\gamma_1 *_a
\gamma_2)(t)=\gamma_1(\frac{t}{a})$
and if $a\leq t\leq 1$, $(\gamma_1 *_a
\gamma_2)(t)=\gamma_2(\frac{t-a}{1-a})$.
\end{nota}

Let us notice that if $\gamma_1$ and $\gamma_2$ are two
non-constant execution paths of a globular complex $X$, then
$\gamma_1 *_a \gamma_2$ is a non-constant execution path of $X$ as
well for any $0<a<1$.

\begin{nota} If $X$ is a globular complex, let ${\P} X:={\P} cat(X)$. \end{nota}

\bth The functor $cat:\gltop\longrightarrow \dtop$ induces a
natural transformation $p:{\P}^{ex}\longrightarrow {\P}$ characterized
by the following facts:
\begin{enumerate}
\item if $X=\glob^{top}(Z)$, then $p_{\glob^{top}(Z)}(t\mapsto (z,t))=z$ for
any $z\in Z$
\item if $\phi\in \gltop(\vI^{top},\vI^{top})$, if $\gamma$ is a non-constant
execution path of a globular complex $X$, then $p_X(\gamma\circ \phi)=p_X(\gamma)$
\item if $\gamma_1$ and $\gamma_2$ are two non-constant
execution paths of a globular complex $X$, then $p_X(\gamma_1 *_a
\gamma_2)=p_X(\gamma_1)*p_X(\gamma_2)$ for any $0<a<1$.
\end{enumerate}
\eth

\bpf See Theorem~\ref{fonc}. \epf

\subsection{Quasi-flow}\label{quasi}

In order to write down in a rigorous way the construction of the
functor $cat$ , the notion of \textit{quasi-flow} seems to be
required.

\bd 
A {\rm quasi-flow} $X$ is a set $X^0$ (the $0$-skeleton) together with
a topological space ${\P}^{top}_{\alpha,\beta}X$ (which can be empty)
for any $(\alpha,\beta)\in X^0\p X^0$ and for any
$\alpha,\beta,\gamma\in X^0\p X^0\p X^0$ a continuous map $]0,1[\p
{\P}^{top}_{\alpha,\beta}X \p {\P}^{top}_{\beta,\gamma}X\to
{\P}^{top}_{\alpha,\gamma}X$ sending $(t,x,y)$ to $x*_t y$ and
satisfying the following condition: if $ab=c$ and $(1-c)(1-d)=(1-b)$,
then $(x*_a y) *_b z= x *_c (y *_d z)$ for any $(x,y,z)\in
{\P}^{top}_{\alpha,\beta}X\p {\P}^{top}_{\beta,\gamma}X\p
{\P}^{top}_{\gamma,\delta}X$. A morphism of quasi-flows
$f:X\longrightarrow Y$ is a set map $f^0:X^0\longrightarrow Y^0$
together with for any $(\alpha,\beta)\in X^0\p X^0$, a continuous map
${\P}^{top}_{\alpha,\beta}X\to {\P}^{top}_{f^0(\alpha),g^0(\beta)}Y$
such that $f(x *_t y)=f(x) *_t f(y)$ for any $x,y$ and any $t\in
]0,1[$. The corresponding category is denoted by $\wdtop$. 
\ed

\bth\label{ssc}
\cite{MR96g:18001a,MR1712872} (Freyd's Adjoint Functor Theorem)
Let $A$ and $X$ be locally small categories. Assume that $A$ is
complete. Then a functor $G:A\longrightarrow X$ has a left adjoint if
and only if it preserves all limits and satisfies the following
``Solution Set Condition''. For each object $x\in X$, there is a set
of arrows $f_i:x\longrightarrow G a_i$ such that for every arrow
$h:x\longrightarrow G a$ can be written as a composite $h=Gt\circ f_i$
for some $i$ and some $t:a_i\longrightarrow a$. \eth

\bth The category of quasi-flows is complete
and cocomplete. \eth

\bpf Let $X:I\longrightarrow \wdtop$ be a diagram of
quasi-flows. Then  the limit of this diagram is
constructed as follows:
\begin{enumerate}
\item the $0$-skeleton is $\limproj X^0$
\item let $\alpha$ and $\beta$ be two elements of $\limproj X^0$ and
let $\alpha_i$ and $\beta_i$ be their image by the canonical
continuous map $\limproj X^0\longrightarrow X(i)^0$
\item let ${\P}_{\alpha,\beta}^{top}(\limproj X):=\limproj_i {\P}_{\alpha_i,\beta_i}^{top}X(i)$.
\end{enumerate}
So all axioms required for the family of topological spaces
${\P}_{\alpha,\beta}^{top}(\limproj X)$ are clearly satisfied. 
Hence the completeness.

The constant diagram functor $\Delta_I$ from
the category of quasi-flows $\wdtop$ to the category
of diagrams of quasi-flows $\wdtop^I$ over a small
category $I$ commutes with limits. It then suffices to
find a set of solutions to prove the existence of a left adjoint by
Theorem~\ref{ssc}. Let $D$ be an object of $\wdtop^I$ and let
$f:D\longrightarrow \Delta_I Y$ be a morphism in $\wdtop^I$. Then one can
suppose that the cardinal $\card(Y)$ of the underlying topological space
$Y^0\sqcup(\bigsqcup_{(\alpha,\beta)\in X^0,X^0}
{\P}^{top}_{\alpha,\beta}Y)$ of $Y$ is lower than the cardinal
$M:=\sum_{i\in I}\card(D(i))$ where $\card(D(i))$ is the cardinal of the
underlying topological space of the quasi-flow $D(i)$. Then let
$\{Z_i,i\in I\}$ be the set of isomorphism classes of
quasi-flows whose underlying topological space is of
cardinal lower than $M$. Then to describe $\{Z_i,i\in I\}$, one has to
choose a $0$-skeleton among $2^M$ possibilities,
for each pair $(\alpha,\beta)$ of the
$0$-skeleton, one has to choose a topological
space among $2^M\p 2^{(2^M)}$ possibilities, and maps $*_t$ among
$(2^{(M\p M\p M)})^{(2^{\aleph_0})}$ possibilities.  Therefore the
cardinal $\card(I)$ of $I$ satisfies
\[\card(I)\leq 2^M \p M\p M \p 2^M\p 2^{(2^M)} \p (2^{(M\p M\p M)})^{(2^{\aleph_0})}\]
so the class $I$ is actually a set. Therefore the class $\bigcup_{i\in
I}\wdtop(D,\Delta_I(Z_i))$ is a set as well. \epf

There is a canonical embedding functor from the category of flows to
that of quasi-flows by setting $*_t=*$ (the composition law of the
flow).

\subsection{Associating a quasi-flow with any globular complex}

\bp\label{key} Let $M$ be a topological space. Let $\gamma_1$ and
$\gamma_2$ be two continuous maps from $[0,1]$ to $M$ with
$\gamma_2(1)=\gamma_1(0)$. Let
$\gamma_3:[0,1]\longrightarrow M$ be another continuous map with
$\gamma_2(1)=\gamma_3(0)$. Assume that
$a,b,c,d\in ]0,1[$ such that $ab=c$ and
$(1-c)(1-d)=(1-b)$. Then
$(\gamma_1 *_a \gamma_2) *_b \gamma_3=\gamma_1 *_c
(\gamma_2 *_d \gamma_3)$. \ep

\bpf Let us calculate $((\gamma_1 *_a \gamma_2) *_b \gamma_3)(t)$.
There are three possibilities:
\begin{enumerate}
\item $0\leq t\leq ab$. Then $((\gamma_1 *_a \gamma_2) *_b
\gamma_3)(t)=\gamma_1(\frac{t}{ab})$.
\item $ab\leq t\leq b$. Then $((\gamma_1 *_a \gamma_2) *_b
\gamma_3)(t)=\gamma_2(\frac{\frac{t}{b}-a}{1-a})=\gamma_2(\frac{t-ab}{b(1-a)})$
\item $b\leq t\leq 1$. Then $((\gamma_1 *_a \gamma_2) *_b
\gamma_3)(t)=\gamma_3(\frac{t-b}{1-b})$
\end{enumerate}

Let us now calculate $(\gamma_1 *_c (\gamma_2 *_d \gamma_3))(t)$.
There are again three possibilities:
\begin{enumerate}
\item $0\leq t\leq c$. Then $(\gamma_1 *_c (\gamma_2 *_d
\gamma_3))(t)=\gamma_1(\frac{t}{c})$
\item $0\leq \frac{t-c}{1-c}\leq d$, or equivalently $c\leq t\leq
c+d(1-c)$. Then $(\gamma_1 *_c (\gamma_2 *_d
\gamma_3))(t)=\gamma_2(\frac{t-c}{d(1-c)})$
\item $d\leq \frac{t-c}{1-c}\leq 1$, or equivalently $c+d(1-c)\leq
t\leq 1$. Then $(\gamma_1 *_c (\gamma_2 *_d
\gamma_3))(t)=\gamma_3(\frac{\frac{t-c}{1-c}-d}{1-d})
=\gamma_3(\frac{t-c-d(1-c)}{(1-d)(1-c)})$.
\end{enumerate}

From $(1-c)(1-d)=(1-b)$, one deduces that $1-c-(1-b)=d(1-c)$, so
$d(1-c)=b-c=b-ab=b(1-a)$. Therefore $d(1-c)=b(1-a)$. So
$c+d(1-c)=b$. The last two equalities complete the proof. \epf

\bp Let $X$ be a globular complex. Let $qcat(X):=X^0$ and
\[{\P}^{top}_{\alpha,\beta}qcat(X):={\P}^{ex}_{\alpha,\beta}X.\] for
any $(\alpha,\beta)\in X^0\p X^0$. This defines a functor
$qcat:\gltop\longrightarrow \wdtop$. \ep

\bpf Immediate consequence of Proposition~\ref{key}. \epf

\bp\label{preserveex} Let $X$ be a globular complex with globular decomposition
\[(\de Z_\beta, Z_\beta,\phi_\beta)_{\beta<\lambda}.\] Then:
\begin{enumerate}
\item for any $\beta<\lambda$, one has the pushout of quasi-flows
\[\xymatrix{qcat(\glob^{top}(\de Z_\beta))\fr{qcat(\phi_\beta)}\fd{}& qcat(X_\beta)\fd{}\\
qcat(\glob^{top}( Z_\beta))\fr{} &
qcat(X_{\beta+1})\cocartesien}\]
\item for any limit ordinal $\beta<\lambda$, the canonical morphism of quasi-flows
\[\liminj_{\alpha<\beta} qcat(X_\alpha)\longrightarrow
qcat(X_\beta)\] is an isomorphism of quasi-flows
\end{enumerate}
\ep

\bpf The first part is a consequence of Proposition~\ref{key}. For any
globular complex $X$, the continuous map $|X_\beta|\longrightarrow
|X_{\beta+1}|$ is a Hurewicz cofibration, and in particular a closed
inclusion of topological spaces. Since $[0,1]$ is compact, it is
$\aleph_0$-small relative to closed inclusions of topological spaces
\cite{MR99h:55031}. Since $\beta$ is a limit ordinal, then $\beta\geq
\aleph_0$. Therefore any continuous map $[0,1]\longrightarrow X_\beta$
factors as a composite $[0,1]\longrightarrow X_\alpha\longrightarrow
X_\beta$ for some $\alpha<\beta$. Hence the second part of the
statement. \epf

\subsection{Construction of the functor $cat$ on objects}

Let $X$ be a globular complex with globular decomposition $(\de
 Z_\beta, Z_\beta,\phi_\beta)_{\beta<\lambda}$. We are going to
construct by induction on $\beta$ a flow $cat(X_\beta)$ and a
morphism of quasi-flows $p_{X_\beta}:qcat(X_\beta)\longrightarrow
cat(X_\beta)$.

There is nothing to do if $X=X_0=(X^0,X^0)$ is a discrete
globular complex. If $X=\glob^{top}(Z)$, then $qcat(X)^0=\{0,1\}$
and \[{\P}^{top}_{0,1}qcat(X)=Z\p \tgltop(\vI^{top},\vI^{top})\]
by Proposition~\ref{decomp}. The projection
${\P}^{top}_{0,1}qcat(X)\longrightarrow Z$ yields a morphism of
quasi-flows $p_X:qcat(X)\longrightarrow cat(X)$.

Let us consider the pushout of multipointed spaces
\[\xymatrix{
\glob^{top}(\de Z_\beta)\fr{\phi_\beta}\fd{} & \fd{}X_\beta\\
\glob^{top}( Z_\beta)\fr{}& \cocartesien X_{\beta+1}}
\]

Let us suppose $p_{X_\beta}:qcat(X_\beta)\longrightarrow cat(X_\beta)$
constructed. Let us consider the set map $i_Z:Z\longrightarrow
{\P}^{ex}\glob^{top}(Z)$ defined by $i_Z(z)(t)=(z,t)$. It is
continuous since it corresponds, by the set map
$\top(Z,{\P}^{ex}\glob^{top}(Z))\longrightarrow \top(Z\p
[0,1],|\glob^{top}(Z)|)$, to the continuous map $(z,t)\mapsto (z,t)$.
The composite
\[\xymatrix@1{ \de Z_\beta\fr{i_{\de Z_\beta}}&{\P}^{ex}\glob^{top}(\de
Z_\beta)\fr{}&
{\P}^{top}(X_\beta)\ar@{->}[rr]^-{qcat(\phi_\beta)}&& {\P}
X_\beta}\] yields a morphism of flows $cat(\phi_\beta):\glob(\de
Z_\beta)\longrightarrow cat(X_\beta)$. Then let
$cat(X_{\beta+1})$ be the flow defined by the pushout of flows
\[\xymatrix{
\glob(\de Z_\beta)\fr{\phi_\beta}\fd{cat(\phi_\beta} & \fd{}cat(X_\beta)\\
\glob( Z_\beta)\fr{}& \cocartesien cat(X_{\beta+1})}
\]
The morphisms of quasi-flows \[qcat(X_\beta)\longrightarrow
cat(X_\beta)\] and
\[qcat(\glob^{top}(Z_\beta))\longrightarrow
cat(\glob^{top}(Z_\beta))\] induce a commutative square of
quasi-flows
\[\xymatrix{qcat(\glob^{top}(\de Z_\beta))\fr{}\fd{}& cat(X_\beta)\fd{}\\
qcat(\glob^{top}( Z_\beta))\fr{} & cat(X_{\beta+1})}\] and
therefore a morphism of quasi-flows
$p_{X_{\beta+1}}:qcat(X_{\beta+1})\longrightarrow
cat(X_{\beta+1})$. If $\beta$ is a limit ordinal, then
$cat(X_\alpha)$ and the morphism of flows
$p_{X_\alpha}:qcat(X_\alpha)\longrightarrow cat(X_\alpha)$ are
defined by induction hypothesis for any $\alpha<\beta$. Then let
$cat(X_\beta):=\liminj_{\alpha<\beta}cat(X_\alpha)$ and
$p_{X_\beta}:=\liminj_{\alpha<\beta}p_{X_\alpha}$.

\subsection{Construction of the functor $cat$ on arrows}

Let $f:X\longrightarrow U$ be a morphism of globular complexes.
The purpose of this section is the construction of
$cat(f):cat(X)\longrightarrow cat(U)$.

If $X=X^0$, then there is nothing to do since the set map
$\gltop(X,U)\longrightarrow \dtop(X,U)$ is just the identity of
$\set(X^0,U^0)$.

If $X=\glob^{top}(Z)$ for some compact space $Z$, let
$f:\glob^{top}(Z)\longrightarrow U$ be a morphism of globular
complexes. Let $cat(f)=p_U\circ qcat(f)\circ i_Z$. Then the mapping
$f\mapsto cat(f)$ yields  a set map
$\gltop(\glob^{top}(Z),U)\longrightarrow
\dtop(\glob(Z),cat(U))$.

Take now a general globular complex $X$ with globular decomposition
\[(\de Z_\beta, Z_\beta,\phi_\beta)_{\beta<\lambda}.\] Using
Theorem~\ref{limittop} and Theorem~\ref{commute}, one obtains a set
map \[\gltop(X_\beta,U)\longrightarrow
\dtop(cat(X_\beta),cat(U))\] and by passage to the limit, a set
map \[cat:\gltop(X,U)\longrightarrow \dtop(cat(X),cat(U)).\]

\subsection{Functoriality of the functor $cat$}

\bth\label{fonc} The mapping $cat(-)$ becomes a functor from
$\gltop$ to $\dtop$. The mapping $p_X:qcat(X)\longrightarrow cat(X)$ yields a
natural transformation $p:qcat\longrightarrow cat$. The mapping
$p_X:{\P}^{top}X\longrightarrow {\P} X$ yields a natural transformation
$p:{\P}^{top}\longrightarrow {\P}$. \eth

\bpf Let $U$ and $V$ be two topological spaces. 
Let $h:U\longrightarrow V$ be a continuous map. 
Let $Z$ be a topological space. Then the following
diagram is clearly commutative:
\[
\xymatrix{
\gltop(\glob^{top}(Z),\glob^{top}(U))\fr{cat(-)}\fd{} & \dtop(\glob(Z),\glob(U))\fd{}\\
\gltop(\glob^{top}(Z),\glob^{top}(V)) \fr{cat(-)} &
\dtop(\glob(Z),\glob(V)) }
\]
where the horizontal maps are both defined by the above construction
and where the right vertical map
$\dtop(\glob(Z),\glob(U))\longrightarrow \dtop(\glob(Z),\glob(V))$ is
induced by the composition by $\glob(h)$.

So for any morphism  $h:U\longrightarrow V$ of globular complexes and
for any topological space $Z$, one has the following commutative diagram
\[
\xymatrix{
\gltop(\glob^{top}(Z),U)\fr{cat(-)}\fd{} & \dtop(\glob(Z),cat(U))\fd{}\\
\gltop(\glob^{top}(Z),V) \fr{cat(-)} & \dtop(\glob(Z),cat(V)) }
\]
where both horizontal maps are defined by the above construction
and where the right vertical map
$\dtop(\glob(Z),cat(U))\longrightarrow \dtop(\glob(Z),cat(V))$ is
induced by the composition by $cat(h)\in
\dtop(cat(U),cat(V))\iso \limproj \dtop(cat(U_\beta),cat(V))$.  Indeed
locally, we are reduced to the situation of the first square.

Take now a general globular complex $X$ with globular
decomposition \[(\de Z_\beta,
Z_\beta,\phi_\beta)_{\beta<\lambda}.\] Then using
Theorem~\ref{limittop} and Theorem~\ref{commute}, one immediately proves
by transfinite induction on $\beta$ that the diagram
\[
\xymatrix{
 \gltop(X_\beta,U)\fr{cat(-)}\fd{} & \dtop(cat(X_\beta),cat(U))\fd{}\\
\gltop(X_\beta,V) \fr{cat(-)} &
\dtop(cat(X_\beta),cat(V)) }
\]
is commutative for any ordinal $\beta<\lambda$. So one obtains the
following commutative diagram
\[
\xymatrix{
\gltop(X,U)\fr{cat(-)}\fd{} & \dtop(cat(X),cat(U))\fd{}\\
\gltop(X,V) \fr{cat(-)} & \dtop(cat(X),cat(V)) }
\]
where both horizontal maps are defined by the above construction
and where the right vertical map
$\dtop(\glob(Z),cat(U))\longrightarrow \dtop(\glob(Z),cat(V))$ is
induced by the composition by $cat(h)\in
\dtop(cat(U),cat(V))\iso \limproj \dtop(cat(U_\beta),cat(V))$.  This
is exactly the functoriality of $cat(-)$.

By specializing the second square to $Z=\{*\}$ and by
Proposition~\ref{pathtopS} and Proposition~\ref{can}, one obtains the
commutative square of topological spaces
\[
\xymatrix{ {\P}^{top}U \fd{{\P}^{top}h} \fr{p_U} & {\P} U \fd{{\P} h}\\
{\P}^{top}V \fr{p_V} & {\P} V}
\]
\epf

\section{Pushout of $\glob(\de Z)\longrightarrow \glob(Z)$ in $\dtop$}\label{expush}

Let $\de Z\longrightarrow Z$ be a continuous map.
Let us consider a diagram of flows as follows:
\[
\xymatrix{
\glob(\de Z) \fr{\phi}\fd{} & A  \fd{}\\
\glob(Z) \fr{} & \cocartesien X}
\]
This short section is devoted to an explicit description of the
pushout $X$ in the category of flows.

Let us consider the set $\mathcal{M}$ of finite sequences
$\alpha_0\dots\alpha_p$ of elements of $A^0=X^0$ with $p\geq 1$
and such that, for any $i$ with $0\leq i\leq p-2$, at least one of
the two pairs $(\alpha_i,\alpha_{i+1})$ and
$(\alpha_{i+1},\alpha_{i+2})$ is equal to $(\phi(0),\phi(1))$.
Let us consider the pushout diagram of topological spaces
\[
\xymatrix{
\de Z \fr{\phi}\fd{} & {\P}_{\phi(0),\phi(1)} A \fd{} \\
Z \fr{} & \cocartesien T }
\]

Let $Z_{\alpha,\beta}={\P}_{\alpha,\beta}A$ if $(\alpha,\beta)\neq
(\phi(0),\phi(1))$ and let $Z_{\phi(0),\phi(1)}=T$. At last, for
any $\alpha_0\dots\alpha_p\in \mathcal{M}$, let
$[\alpha_0\dots\alpha_p]= Z_{\alpha_0,\alpha_1}\p
Z_{\alpha_1,\alpha_2}\p \dots \p Z_{\alpha_{p-1},\alpha_p}$.  And
$[\alpha_0\dots\alpha_p]_i$ denotes the same product as
$[\alpha_0\dots\alpha_p]$ except that
$(\alpha_i,\alpha_{i+1})=(\phi(0),\phi(1))$ and that the factor
$Z_{\alpha_{i},\alpha_{i+1}}=T$ is replaced by
${\P}_{\phi(0),\phi(1)} A$.  That means that in the product
$[\alpha_0\dots\alpha_p]_i$, the factor ${\P}_{\phi(0),\phi(1)} A$
appears exactly once. For instance, one has (with $\phi(0)\neq
\phi(1)$) \beas &&
[\alpha\phi(0)\phi(1)\phi(0)\phi(1)]={\P}_{\alpha,\phi(0)}A\p T\p
{\P}_{\phi(1),\phi(0)}A\p T\\
&& [\alpha\phi(0)\phi(1)\phi(0)\phi(1)]_1={\P}_{\alpha,\phi(0)}A\p
{\P}_{\phi(0),\phi(1)} A\p {\P}_{\phi(1),\phi(0)}A\p T\\
&&  [\alpha\phi(0)\phi(1)\phi(0)\phi(1)]_3={\P}_{\alpha,\phi(0)}A\p
T\p {\P}_{\phi(1),\phi(0)}A\p {\P}_{\phi(0),\phi(1)} A. \eeas The
idea is that in the products $[\alpha_0\dots\alpha_p]$, there are
no possible simplifications using the composition law of $A$. On
the contrary, exactly one simplification is possible using the
composition law of $A$ in the products
$[\alpha_0\dots\alpha_p]_i$. For instance, with the examples
above, there exist continuous maps
\[[\alpha\phi(0)\phi(1)\phi(0)\phi(1)]_1\longrightarrow
[\alpha\phi(0)\phi(1)]\] and
\[[\alpha\phi(0)\phi(1)\phi(0)\phi(1)]_3\longrightarrow[\alpha\phi(0)\phi(1)\phi(1)]\]
induced by the composition law of $A$ and there exist
continuous maps
\[[\alpha\phi(0)\phi(1)\phi(0)\phi(1)]_1\longrightarrow [\alpha\phi(0)\phi(1)\phi(0)\phi(1)]\]
and
\[[\alpha\phi(0)\phi(1)\phi(0)\phi(1)]_3\longrightarrow [\alpha\phi(0)\phi(1)\phi(0)\phi(1)]\]
induced by the continuous map ${\P}_{\phi(0),\phi(1)}A\longrightarrow T$.

Let ${\P}_{\alpha,\beta}M$ be the colimit of the diagram of
topological spaces consisting of the topological spaces
$[\alpha_0\dots\alpha_p]$ and $[\alpha_0\dots\alpha_p]_i$ with
$\alpha_0=\alpha$ and $\alpha_p=\beta$ and with the two kinds of maps
above defined. The composition law of $A$ and the free
concatenation obviously gives  a continuous associative map
${\P}_{\alpha,\beta}M\p {\P}_{\beta,\gamma}M\longrightarrow
{\P}_{\alpha,\gamma}M$.

\bp \label{pushexplicit}\label{pushexplicit0} (\cite{model3} Proposition~15.1)
One has the pushout diagram of flows
\[
\xymatrix{
\glob(\de Z) \fr{\phi}\fd{} & A \fd{} \\
\glob(Z) \fr{} & \cocartesien M }
\]
\ep

\section{Geometric realization of execution paths}\label{stretched}

\bp \label{span} Let $Z$ be a compact topological space. Let $f$
and $g$ be two morphisms of globular complexes from $\glob^{top}(Z)$
to a globular complex $U$ such that the continuous maps ${\P} f$ and
${\P} g$ from $Z$ to ${\P} U$ are equal. Then there exists one and
only one map $\phi:|\glob^{top}(Z)|\longrightarrow [0,1]$ such that
\[f((z,t))= {\P}^{top} g(t\mapsto
(z,t))(\phi(z,t)).\] 
Moreover this map $\phi$ is necessarily continuous. \ep

Notice that the map $\phi:\glob^{top}(Z)\longrightarrow [0,1]$
induces a morphism of globular complexes from
$\glob^{top}(Z)$ to $\vI^{top}$.

\bpf By hypothesis, the equality ${\P}^{top}
f([0,1])={\P}^{top} g([0,1])$ holds. For a
given $z_0\in Z$, if
\[\{0=t_0<\dots<t_p=1\}={\P}^{top} f(t\mapsto (z_0,t))([0,1])\cap
U^0\] and
\[\{0=t'_0<\dots<t'_p=1\}={\P}^{top} g(t\mapsto (z_0,t))([0,1])\cap U^0\] then necessarily
$\phi(z_0,t_i)=t'_i$ for  $0\leq i \leq p$. For $t\in
]t_i,t_{i+1}[$, the map
\[{\P}^{top} g(t\mapsto (z_0,t))\!\restriction_{]t'_i,t'_{i+1}[}\]
is one-to-one by hypothesis. Therefore for $t\in ]t_i,t_{i+1}[$,
${\P}^{top} f(t\mapsto
(z_0,t))(t)$ is equal to
{
\[  ( {\P}^{top} g(t\mapsto (z_0,t))\!\restriction_{]t'_i,t'_{i+1}[})( {\P}^{top} g(t\mapsto (z_0,t))\!\restriction_{]t'_i,t'_{i+1}[})^{-1}{\P}^{top} f(t\mapsto (z_0,t))(t)\]}
so necessarily one has
\[\phi(z_0,t)=( {\P}^{top} g(t\mapsto (z_0,t))\!\restriction_{]t'_i,t'_{i+1}[})^{-1}{\P}^{top} f(t\mapsto (z_0,t))(t)\]
Now suppose that $\phi$ is not continuous at
$(z_\infty,t_\infty)$. Then there exists an open
neighborhood $U$ of $\phi(z_\infty,t_\infty)$ such
that for any open $V$ containing
$(z_\infty,t_\infty)$, for any $(z,t)\in
V\backslash\{(z_\infty,t_\infty)\}$,
$\phi(z,t)\notin U$. Take a sequence
$(z_n,t_n)_{n\geq 0}$ of $V$ tending to
$(z_\infty,t_\infty)$. Then there exists a subsequence
of $(\phi(z_n,t_n))_{n\geq 0}$ tending
to some $t'\in [0,1]$ since $[0,1]$ is compact: by hypothesis
$t'$ is in the topological closure of the complement of $U$; this
latter being closed, $t'\notin U$. So we can take
$(z_n,t_n)_{n\geq 0}$ such that
$(\phi(z_n,t_n))_{n\geq 0}$ converges. Then
$f((z_n,t_n))$ tends to
$f((z_\infty,t_\infty))$ because $f$ is
continuous, ${\P}^{top} g(t\mapsto (z_n,t))$
tends to ${\P}^{top} g(t\mapsto
(z_\infty,t))$ for the Kelleyfication of the
compact-open topology so
$f((z_\infty,t_\infty))= {\P}^{top}
g(t\mapsto (z_\infty,t))(t')$
with $t' \notin U$ and $t'=\phi(z_\infty,t_\infty)\in
U$: contradiction. \epf

\bth\label{i} For any globular complex $X$, there exists a
continuous map $i_X:{\P} X\longrightarrow {\P}^{top} X$ such that
$p_X\circ i_X=\id_{{\P} X}$. \eth

Notice that $i_X$ cannot be obtained from a morphism of
quasi-flows. Otherwise one would have $(x
*_a y) *_a z=x *_a (y *_a z)$ in ${\P}^{top}X$ for
some fixed $a\in ]0,1[$, and this is impossible.

\bpf First of all, notice that there is an inclusion of sets
$\top({\P} X,{\P}^{top}X)\subset \top({\P} X\p[0,1],X)$. So
constructing a continuous map from ${\P} X$ to ${\P}^{top}X$ is
equivalent to constructing a continuous map from ${\P} X\p[0,1]$ to
$X$ satisfying some obvious properties, since the category $\top$ of
compactly generated topological spaces is cartesian closed.

Let $X$ be a globular complex with globular decomposition $(\de
 Z_\beta, Z_\beta,\phi_\beta)_{\beta<\lambda}$. We are going to
construct a continuous map $i_{X_\beta}:{\P}
X_\beta\longrightarrow {\P}^{top} X_\beta$. For $\beta=0$, there
is nothing to do since the topological spaces are both discrete.
Assume that $i_{X_\beta}:{\P} X_\beta\longrightarrow {\P}^{top}
X_\beta$ is constructed for some $\beta\geq 0$ such that
$p_{X_\beta}\circ i_{X_\beta}=\id_{{\P} X_\beta}$. Let us
consider the pushout of multipointed spaces
\[\xymatrix{
\glob^{top}(\de Z_\beta)\fr{\phi_\beta}\fd{} & \fd{}X_\beta\\
\glob^{top}( Z_\beta)\fr{\overline{\phi_\beta}}& \cocartesien X_{\beta+1}}
\]

Proposition~\ref{pushexplicit0} provides an explicit method for the
calculation of ${\P} X_{\beta+1}$ as the colimit of a diagram of
topological spaces. Let us consider the pushout diagram
of topological spaces
\[
\xymatrix{
\de Z_\beta \fr{\phi_\beta}\fd{} & {\P}_{\phi_\beta(0),\phi_\beta(1)} X_\beta \fd{} \\
Z_\beta \fr{} & \cocartesien T }
\]
Constructing a continuous map ${\P}
X_{\beta+1}\longrightarrow {\P}^{top} X_{\beta+1}$ is then equivalent to
constructing continuous maps $[\alpha_0\dots\alpha_p]\longrightarrow
{\P}^{top} X_{\beta+1}$ and
$[\alpha_0\dots\alpha_p]_i\longrightarrow
{\P}^{top} X_{\beta+1}$ for any finite sequence
$\alpha_0\dots\alpha_p$ of $\mathcal{M}$ such that any diagram
like
\[
\xymatrix{
[\alpha_0\dots\alpha_p]_i\fd{} \fr{}&{\P}^{top} X_{\beta+1}&&[\alpha_0\dots\alpha_p]_i\fd{} \fr{}
&{\P}^{top} X_{\beta+1}\\
[\alpha_0\dots\alpha_p]\ar@{->}[ru]&&&[\alpha_0\dots \widehat{\phi(0)\phi(1)}\dots\alpha_p]\ar@{->}[ru]&
}
\]
is commutative.

We are going to proceed by induction on $p$. If $p=1$, then
$[\alpha_0\alpha_1]$ is equal to ${\P}_{\alpha_0,\alpha_1}X_\beta$ if
$(\alpha_0,\alpha_1)\neq (\phi_\beta(0),\phi_\beta(1))$ and is equal
to $T$ if $(\alpha_0,\alpha_1)= (\phi_\beta(0),\phi_\beta(1))$. For
$p=1$, the only thing we then have to prove is that the continuous map
$p_{X_{\beta+1}}:{\P}^{top}X_{\beta+1}\longrightarrow {\P} X_{\beta+1}$
has the right lifting property with respect to the continuous map
${\P}_{\phi_\beta(0),\phi_\beta(1)}X_\beta\longrightarrow T$, in other terms 
that there exists a continuous map $k:T\longrightarrow {\P}^{top}X_{\beta+1}$ making
commutative the diagram of topological spaces
\[
\xymatrix{
{\P}_{\phi_\beta(0),\phi_\beta(1)}X_\beta \fd{} \fr{i_{X_\beta}} & {\P}^{top}X_{\beta+1}\ar@{->}[d]^-{p_{X_{\beta+1}}}\\
T \fr{} \ar@{-->}[ru]^-{k} & {\P} X_{\beta+1}}
\]
Since ${\P}_{\phi_\beta(0),\phi_\beta(1)}X_\beta\longrightarrow T$ is a
pushout of a NDR pair of spaces, then the pair of spaces
$(T,{\P}_{\phi_\beta(0),\phi_\beta(1)} X_\beta)$ is a NDR pair as
well. If $z\in Z_\beta$, let $[z](t)=(z,t)$ for $t\in [0,1]$. This defines
an execution path of $\glob^{top}(Z_\beta)$. Then
$\overline{\phi_\beta}\circ [z]$ is still an execution path. Since
${\P}_{\phi_\beta(0),\phi_\beta(1)}X_\beta\longrightarrow T$ is a
(closed) inclusion of topological spaces, then for any $z\in \de
Z_\beta$, $\overline{\phi_\beta}\circ [z]$ is an execution
path of $X_\beta$. By Proposition~\ref{span} and since $\de Z_\beta$
is compact, there exists a continuous map $\psi:\de
Z_\beta\p[0,1]\longrightarrow [0,1]$ such that
\[i_{X_\beta}(z)(t)=(\overline{\phi_\beta}\circ [z])(\psi(z,t)).\]
Then define $k$ by: $k(x)=i_{X_\beta}(x)$ if $x\in
{\P}_{\psi(0),\psi(1)}X_\beta$ and
\[k(x)(t)=(\overline{\phi_\beta}\circ
[x])(\mu(x)t+(1-\mu(x))\psi(x,t))\] if $x\in Z_\beta\backslash
\de  Z_\beta$. The case $p=1$ is complete.

We now have to construct $[\alpha_0\dots\alpha_p]_i\longrightarrow
{\P}^{top}X_{\beta+1}$ and $[\alpha_0\dots\alpha_p]\longrightarrow
{\P}^{top}X_{\beta+1}$ by induction on $p\geq 1$. The product
$[\alpha_0\dots \widehat{\phi(0)\phi(1)}\dots\alpha_p]$ is of length
strictly lower than $p$. Therefore the continuous map $[\alpha_0\dots
\widehat{\phi(0)\phi(1)}\dots\alpha_p]\longrightarrow
{\P}^{top}X_{\beta+1}$ is already constructed. Then the commutativity of
the diagram
\[
\xymatrix{
[\alpha_0\dots\alpha_p]_i\fd{} \fr{}
&{\P}^{top} X_{\beta+1}\\
[\alpha_0\dots \widehat{\phi(0)\phi(1)}\dots\alpha_p]\ar@{->}[ru]&
}
\]
entails the definition of $[\alpha_0\dots\alpha_p]_i\longrightarrow
{\P}^{top}X_{\beta+1}$.
It remains to prove that there exists $k$ making the following diagram
commutative:
\[
\xymatrix{
[\alpha_0\dots\alpha_p]_i \fd{} \fr{} & {\P}^{top}X_{\beta+1}\ar@{->}[d]^-{p_{X_{\beta+1}}}\\
[\alpha_0\dots\alpha_p] \fr{} \ar@{-->}[ru]^-{k} & {\P} X_{\beta+1}}
\]
Once again the closed inclusion
$[\alpha_0\dots\alpha_p]_i\longrightarrow [\alpha_0\dots\alpha_p]$ is
a Hurewicz cofibration. There are three
mutually exclusive possible cases:
\begin{enumerate}
\item $[\alpha_0\dots\alpha_p]_i=P\p {\P}_{\phi_\beta(0),\phi_\beta(1)}X_\beta
\p Q$ and $[\alpha_0\dots\alpha_p]=P\p T \p Q$ where $P$ and $Q$ are
objects of the diagram of topological spaces calculating ${\P}
X_{\beta+1}$.
\item $[\alpha_0\dots\alpha_p]_i=P\p {\P}_{\phi_\beta(0),\phi_\beta(1)}X_\beta$
and $[\alpha_0\dots\alpha_p]=P\p T$ where $P$ is an
object of the diagram of topological spaces calculating ${\P}
X_{\beta+1}$.
\item $[\alpha_0\dots\alpha_p]_i={\P}_{\phi_\beta(0),\phi_\beta(1)}X_\beta
\p Q$ and $[\alpha_0\dots\alpha_p]=T \p Q$ where $Q$ is
an object of the diagram of topological spaces calculating ${\P}
X_{\beta+1}$.
\end{enumerate}
Let us treat for instance the first case. The products $P$ and $Q$ are
of length strictly lower than $p$. So by induction hypothesis,
$i_{X_{\beta+1}}:P\longrightarrow {\P}^{top}X_{\beta+1}$ and
$i_{X_{\beta+1}}:Q\longrightarrow {\P}^{top}X_{\beta+1}$ are already
constructed. For any $z\in Z_\beta$ and any $(p,q)\in P\p Q$, consider the
execution path
\[\Gamma(p,z,q):=(i_{X_{\beta+1}}(p)*_{1/2} (\overline{\phi_\beta}\circ [z]))*_{1/2} i_{X_{\beta+1}}(q).\]
By Proposition~\ref{span} and since $\de Z_\beta$ is
compact, there exists a continuous map $\psi:\de P\p Z_\beta\p
Q\p[0,1]\longrightarrow [0,1]$ such that
$\Gamma(p,z,q)(\psi(p,z,q,t))=i_{X_\beta}(p,z,q)(t)$. Then define $k$ by:
\begin{enumerate}
\item
$k(p,x,q)=i_{X_\beta}(p,x,q)(t)$ if $x\in
{\P}_{\psi(0),\psi(1)}X_\beta$
\item $k(p,x,q)=\Gamma(p,z,q)(\mu(x)t+(1-\mu(x))\psi(p,x,q,t))$
if $x\in Z_\beta\backslash \de  Z_\beta$.
\end{enumerate}
The induction is complete. \epf

\section{Conclusion}

Since the functor $cat:\gltop\longrightarrow \dtop$ is constructed, we
are now ready to compare the S-homotopy equivalences in the two
frameworks.

\part{S-homotopy and  flow}\label{model2}

\section{Introduction}

Section~\ref{SHEP} studies the notion of \textit{S-homotopy extension
property} for morphisms of globular complexes.  This is the analogue
in our framework of the notion of Hurewicz cofibration.  This section,
as short as possible, studies some analogues of well-known theorems in
homotopy theory of topological spaces.  The goal of
Section~\ref{comparing} is the comparison of the space of morphisms of
globular complexes from a globular complex $X$ to a globular complex
$U$ with the space of morphisms of flows from the flow $cat(X)$ to the
flow $cat(U)$. It turns out that these two spaces are homotopy
equivalent.  The proof requires the careful study of two transfinite
towers of topological spaces and needs the introduction of a model
category of topological spaces which is not the usual one, but another
one whose weak equivalences are the homotopy equivalences
\cite{MR35:2284} \cite{MR39:4846} \cite{ruse}.  At last,
Section~\ref{compS} makes the comparison between the two notions of
S-homotopy equivalences using all previous results.

\section{S-homotopy extension property}\label{SHEP}

We first need to develop some of the theory of morphisms of globular
complexes satisfying the S-homotopy extension property in order to
obtain Corollary~\ref{montage2}.

\bd 
Let $i:A\longrightarrow X$ be a morphism of globular complexes and let
$Y$ be a globular complex. The morphism $i:A\longrightarrow X$
satisfies {\rm the S-homotopy extension property} for $Y$ if for any
morphism $f:X\longrightarrow Y$ and any S-homotopy $h:[0,1]\boxtimes
A\longrightarrow Y$ such that for any $a\in A$, $h(0\boxtimes
a)=f(i(a))$, there exists a S-homotopy $H:[0,1]\boxtimes
X\longrightarrow Y$ such that for any $x\in X$, $H(0\boxtimes x)=f(x)$
and for any $(t,a)\in [0,1]\p A$, $H(t\boxtimes i(a))=h(t\boxtimes
a)$. \ed

\bd A morphism of globular complexes $i:A\longrightarrow X$ satisfies the
{\rm S-homotopy extension property} if $i:A\longrightarrow X$
satisfies the S-homotopy extension property for any globular complex
$Y$.
\ed

\bp
Let $i:A\longrightarrow X$ be a morphism of globular complexes. Let us
consider the cocartesian diagram of multipointed topological spaces
\[
\xymatrix{
\{0\}\boxtimes A \fr{}\fd{i} & [0,1]\boxtimes A\fd{}\\
\{0\}\boxtimes X \fr{} & \cocartesien Mi}
\]
Then $Mi$ inherits a globular decomposition from those of $A$ and
$X$. This makes the multipointed topological space $Mi$ into a
globular complex. Moreover both morphisms $X\longrightarrow Mi$ and
$[0,1]\boxtimes A\longrightarrow Mi$ are morphisms of globular
complexes. One even has $(Mi)_\beta=X$ for some ordinal $\beta$ and
$X\longrightarrow Mi$ is the canonical morphism induced by the
globular decomposition of $Mi$. \ep

\bpf Let $(\de Z_\beta, Z_\beta,\phi_\beta)_{\beta<\lambda}$ be the
globular decomposition of $A$. The morphism of multipointed spaces
$\{0\}\boxtimes A\longrightarrow [0,1]\boxtimes A$ can be viewed as a
composite \[\{0\}\boxtimes A\longrightarrow \{0,1\}\boxtimes
A\longrightarrow [0,1]\boxtimes A.\] The morphism of globular
complexes $\{0\}\boxtimes A\longrightarrow \{0,1\}\boxtimes A$ is the
transfinite composition of pushouts of the morphisms $\glob^{top}(\de
Z_\beta)\longrightarrow \glob^{top}(Z_\beta)$ for $\beta<\lambda$.
The morphism of globular complexes $\{0,1\}\boxtimes A\longrightarrow
[0,1]\boxtimes A$ is the transfinite composition of pushouts of the
$\glob^{top}(Z_\beta\sqcup Z_\beta)\longrightarrow \glob^{top}([0,1]\p
Z_\beta)$.  Therefore the morphism of multipointed spaces
$X\longrightarrow Mi$ is a relative globular complex. So $Mi$ has a
canonical structure of globular complexes and both morphisms
$X\longrightarrow Mi$ and $[0,1]\boxtimes A\longrightarrow Mi$ are
morphisms of globular complexes. \epf

The commutative diagram of globular complexes
\[
\xymatrix{
\{0\}\boxtimes A \fr{}\fd{i} & [0,1]\boxtimes A\fd{}\\
\{0\}\boxtimes X \fr{} & [0,1]\boxtimes X}
\]
gives rise to a morphism of multipointed spaces $\psi(i):Mi\longrightarrow [0,1]\boxtimes X$.
Since by Proposition~\ref{preserveex}, one also has the cocartesian diagram of quasi-flows
\[
\xymatrix{
qcat(\{0\}\boxtimes A) \fr{}\fd{i} & qcat([0,1]\boxtimes A)\fd{}\\
qcat(\{0\}\boxtimes X) \fr{} & \cocartesien qcat(Mi)}
\]
then there exists a morphism of quasi-flows
$qcat(Mi)\longrightarrow qcat([0,1]\boxtimes X)$. Therefore the
morphism of multipointed spaces $\psi(i):Mi\longrightarrow
[0,1]\boxtimes X$ satisfies
\[\psi(i)({\P}^{top}Mi)\subset {\P}^{top}([0,1]\boxtimes X).\] 
So $\psi(i)$ is a morphism of globular complexes.

\bth\label{retract} Let $i:A\longrightarrow X$ be a morphism of
globular complexes. Then the following assertions are equivalent:
\begin{enumerate}
\item the morphism $i$ satisfies the S-homotopy extension property
\item the morphism of globular complexes $\psi(i)$ has a retract $r$,
that is to say there exists a morphism of globular complexes
\[r:[0,1]\boxtimes X\longrightarrow ([0,1]\boxtimes A)\sqcup_{\{0\}\boxtimes A} (\{0\}\boxtimes X)\]
such that $r\circ \psi(i)=\Id_{([0,1]\boxtimes
A)\sqcup_{\{0\}\boxtimes A} (\{0\}\boxtimes X)}$.
\end{enumerate}
\eth

The proof is exactly the same as the one of \cite{model3}
Theorem~9.4. The main point is that the multipointed space $Mi$ is a
globular complex.

\bpf Giving two morphisms of globular complexes $f:X\longrightarrow Y$ and
$h:[0,1]\boxtimes A\longrightarrow Y$ such that $h(0\boxtimes
a)=f(i(a))$ for any $a\in A$ is
equivalent to giving a morphism of globular complexes still denoted by $h$
from $([0,1]\boxtimes A)\sqcup_{\{0\}\boxtimes A}
(\{0\}\boxtimes X)$ to $Y$. The
S-homotopy extension problem for $i$ has then
always a solution if and only for any morphism of globular complexes
$h:([0,1]\boxtimes A)\sqcup_{\{0\}\boxtimes A}
(\{0\}\boxtimes X)\longrightarrow Y$, there exists a
morphism of globular complexes $H:[0,1]\boxtimes X\longrightarrow Y$ such that
$H\circ \psi(i)=h$. Take $Y=([0,1]\boxtimes
A)\sqcup_{\{0\}\boxtimes A} (\{0\}\boxtimes X)$
and let $h$ be the identity map of $Y$. This yields the retract
$r$.  Conversely, let $r$ be a retract of $i$. Then $H:=h\circ r$
is always a solution of the S-homotopy
extension problem.  \epf

\bth 
Let $(Z,\de Z)$ be a NDR pair of compact spaces.  Then the inclusion
of globular complexes $i:\glob^{top}(\de Z)\longrightarrow
\glob^{top}(Z)$ satisfies the S-homotopy extension property. \eth

\bpf Since $(Z,\de Z)$ is a NDR pair, then the closed inclusion
$[0,1]\p \de Z\cup \{0\}\p Z\longrightarrow [0,1]\p Z$ has a retract
$[0,1]\p Z\longrightarrow [0,1]\p \de Z\cup \{0\}\p Z$. Then the
morphism of globular complexes $\glob^{top}([0,1]\p \de Z\cup \{0\}\p
Z)\longrightarrow \glob^{top}([0,1]\p Z)$ has a retract
$\glob^{top}([0,1]\p Z)\longrightarrow \glob^{top}([0,1]\p \de Z\cup
\{0\}\p Z)$. Hence the result by Theorem~\ref{retract}. \epf

\bth\label{limtopt} 
Let $U$ be a compact connected non-empty space.  Let $X$ and $Y$ be
two globular complexes. Then there exists a natural homeomorphism
\[\ttop(U,\tgltop(X,Y))\iso \tgltop(U\boxtimes X,Y).\]
\eth

\bpf We already know by Theorem~\ref{existencegl} that there exists
a natural bijection
\[\top(U,\tgltop(X,Y))\iso \gltop(U\boxtimes X,Y).\]
Let $(\de Z_\beta, Z_\beta,\phi_\beta)_{\beta<\lambda}$ be the globular
decomposition of $X$. We are going to prove that 
\[\ttop(U,\tgltop(X_\beta,Y))\iso \tgltop(U\boxtimes X_\beta,Y).\]
Using the construction of $\boxtimes$ and Theorem~\ref{limittop}, it
suffices to prove the homeomorphism for $X=X_0$ and
$X=\glob^{top}(Z)$.  The space $\tgltop(X_0,Y)$ is the discrete space
of set maps $\set(X^0,Y^0)$ from $X^0$ to $Y^0$. Since $U$ is
connected and non-empty, one has the homeomorphism $\ttop(U,\tgltop(X_0,Y))\iso
\set(X^0,Y^0)$. On the other hand,
$\tgltop(U\boxtimes X_0,Y)\iso \tgltop(X_0,Y)\iso \set(X^0,Y^0)$, hence
the result for $X_0$. At last,
\beas
&& \top(W,\ttop(U,\tgltop(\glob^{top}(Z),Y)))\\
&&\iso \top(W\p U,\tgltop(\glob^{top}(Z),Y))\\
&&\iso \gltop((W\p U)\boxtimes \glob^{top}(Z),Y)\\
&&\iso \gltop(\glob^{top}(W\p U\p Z),Y)
\eeas
and $\top(W,\tgltop(U\boxtimes \glob^{top}(Z),Y))\iso
\top(W,\tgltop(\glob^{top}(U\p Z),Y))$.  It is then easy to see that
both sets $\gltop(\glob^{top}(W\p U\p Z),Y)$ and
$\top(W,\linebreak[4]\tgltop(\glob^{top}(U\p Z),Y))$ can be identified
with the same subset of $\top([0,1]\p W\p U\p Z,Y)$. Hence the result
by Yoneda.
\epf

\bth\label{car} 
A morphism of globular complexes $i:A\longrightarrow X$ satisfies the
S-homotopy extension property if and only if for any globular complex
$Y$, the continuous map $i^*:\tgltop(X,Y)\longrightarrow\tgltop(A,Y)$
is a Hurewicz fibration. \eth

\bpf For any topological space $M$, one has
\[\top([0,1]\p M,\tgltop(A,Y))\iso \top (M,\ttop([0,1],\tgltop(A,Y)))\]
since $\top$ is cartesian closed. One also has \[\top
(M,\ttop([0,1],\tgltop(A,Y)))\iso \top(M,\tgltop([0,1]\boxtimes A,Y))\] by
Theorem~\ref{limtopt}. Considering a  commutative diagram like
\[
\xymatrix{
\{0\}\p M\ar@{^{(}->}[d] \fr{\phi} & \tgltop(X,Y)\fd{i^*}\\
[0,1]\p M \fr{\psi}\ar@{-->}[ru]^-{k} & \tgltop(A,Y)}\]
is then equivalent to considering a commutative diagram of topological spaces
\[\xymatrix{
M \fd{}\fr{} & \tgltop(\{0\}\boxtimes X,Y) \fd{}\\
\tgltop([0,1]\boxtimes A,Y)\fr{} & \tgltop(\{0\}\boxtimes A,Y)}
\]
Since $\{0\}\boxtimes A\longrightarrow [0,1]\boxtimes A$ is a relative
globular complex, using again Theorem~\ref{limittop}, considering such
a commutative diagram is equivalent to considering a continuous map
$M\longrightarrow \tgltop(Mi,Y)$. Finding a continuous map $k$ making
both triangles commutative is equivalent to finding a commutative
diagram of the form
\[\xymatrix{
M \fr{\overline{\phi}} \fd{=}&  \tgltop(Mi,Y) \\
M \ar@{-->}[r]^-{\ell}& \tgltop([0,1]\boxtimes X,Y)\fu{\psi(i)^*}}
\]
If $i:A\longrightarrow X$ satisfies the S-homotopy extension property,
then $\psi(i):Mi\longrightarrow [0,1]\boxtimes X$ has a retract
$r:[0,1]\boxtimes X\longrightarrow Mi$.  Then take
$\ell=\overline{\phi}\circ r$. Conversely, if $\ell$ exists for any $M$
and any $Y$, take $M=\{0\}$ and $Y=Mi$ and
$\overline{\phi}(0)=\id_{Mi}$. Then $\ell(0)$ is a retract of
$\psi(i)$. Therefore $i:A\longrightarrow X$ satisfies the S-homotopy
extension property. \epf

\begin{cor} \label{montage2} 
Let $Z$ be a compact space and let $\de Z\subset Z$ be a compact
subspace such that the canonical inclusion is a NDR pair.  Let $U$ be
a globular complex.  Then the canonical restriction map \[
\tgltop(\glob^{top}(Z),U)\to
\tgltop(\glob^{top}(\de Z),U)\] is a
Hurewicz fibration. \end{cor}

\section{Comparing execution paths of globular complexes and of
flows}\label{comparing}

\subsection{Morphisms of globular complexes and morphisms of flows}

\bp \label{surjectivite} 
Let $Z$ be a compact topological space. Let $U$ be a globular
complex. Consider the set map
\[cat:\gltop(\glob^{top}(Z),U)\longrightarrow
\dtop(\glob(Z),cat(U)).\]
\begin{enumerate}
\item The mapping \[cat:\tgltop(\glob^{top}(Z),U)\longrightarrow \tdtop(\glob(Z),cat(U))\] is
continuous.
\item  There exists a continuous  map \[r:\tdtop(\glob(Z),cat(U))\longrightarrow \tgltop(\glob^{top}(Z),U)\]
such that $cat\circ
r=Id_{\tdtop(\glob(Z),cat(U))}$.
In particular, this means that $cat$ is onto.
\item The map $r\circ cat$ is homotopic to $Id_{\tgltop(\glob^{top}(Z),U)}$. In particular, 
this means that $\tgltop(\glob^{top}(Z),U)$
and $\tdtop(\glob(Z),cat(U))$
are homotopy equivalent.
\end{enumerate}
\ep

\bpf One has
\[\tgltop(\glob^{top}(Z),U)\subset \ttop(Z\p [0,1],U)\iso \ttop(Z,\ttop([0,1],U))\]
therefore
\[\tgltop(\glob^{top}(Z),U)\iso \bigsqcup_{(\alpha,\beta)\in U^0\p U^0} \ttop(Z,{\P}_{\alpha,\beta}^{top}U).\]
On the other hand,
\[\tdtop(\glob(Z),cat(U))\iso
\bigsqcup_{(\alpha,\beta)\in U^0\p U^0}
\ttop(Z,{\P}_{\alpha,\beta}U).\] So the set map \[cat:
\tgltop(\glob^{top}(Z),U)\to
\tdtop(\glob(Z),cat(U))\] is
induced by $p_U$ which is continuous. Hence \[
cat:\tgltop(\glob^{top}(Z),U)\to
\tdtop(\glob(Z),cat(U))\] is
continuous. Choose a map $i_U$ like in Theorem~\ref{i}. Let
\[r(\phi)\in
\tgltop(\glob^{top}(Z),U)\] defined by
$r(\phi)((z,t)):=(i_U\phi(z))(t)$.
Then
\begin{alignat*}{2}
cat(r(\phi))(z) &= p_U\circ {\P}^{top}(r(\phi))(t\mapsto (z,t)) && \hbox{ by definition of $cat(-)$}\\
&=p_U \circ {\P}^{top}(r(\phi)) \circ i_{\glob^{top}(Z)}(z)&& \hbox{ by definition of $i_{\glob^{top}(Z)}$}\\
&= p_U(r(\phi)\circ i_{\glob^{top}(Z)}(z)) && \hbox{ by definition of ${\P}^{top}$}\\
&= p_U(t\mapsto r(\phi)((z,t))) && \hbox{ by definition of $i_{\glob^{top}(Z)}$}\\
&= p_U(i_U \phi(z)) && \hbox{ by definition of $r(\phi)$}\\
&= \phi(z) && \hbox{ since $p_U\circ i_U=\Id$}\\
\end{alignat*}
therefore $cat(r(\phi))=\phi$. So the
second assertion holds. One has
\begin{alignat*}{2}
(r\circ cat (f))(z,t) &= (i_U \circ cat (f)(z))(t) && \hbox{ by definition of $r$}\\
&= (i_U \circ p_U \circ {\P}^{top}(f)(t\mapsto
(z,t)))(t)&& \hbox{ by
definition of $cat$}
\end{alignat*}
Since $(i_U \circ p_U \circ
{\P}^{top}(f)(t\mapsto
(z,t)))$ is an execution path of $U$
by Theorem~\ref{i},
since \[p_U\circ (i_U \circ p_U \circ
{\P}^{top}(f)(t\mapsto
(z,t)))= p_U \circ
{\P}^{top}(f)(t\mapsto (z,t))={\P}
f(z),\] then by Proposition~\ref{span}, there exists a
continuous map $\phi:Z\p [0,1]\longrightarrow [0,1]$ such that
\[f((z,t))={\P}^{top}(f)(t\mapsto (z,t))= (r\circ cat (f))(z,\phi(z,t))\]
Notice that for a given $z\in Z$, the mapping $t\mapsto
\phi(z,t)$ is necessarily non-decreasing. Hence the
third assertion by considering the homotopy
\[H(f,u)((z,t))= (r\circ cat (f))(z,u\phi(z,t)+(1-u)t)\]
\epf

\bp (\cite{model3} Corollary~9.9) \label{montage} 
Let $Z$ be a compact space and let $\de Z\subset Z$ be a compact
subspace such that the canonical inclusion is a NDR pair.  Let $U$ be
a flow. Then the canonical restriction map 
\[
\tdtop(\glob(Z),U)\to
\tdtop(\glob(\de Z),U)\] is a Hurewicz fibration. \ep

\subsection{Homotopy limit of a transfinite tower and homotopy pullback}

Corollary~\ref{tourtransfinie} and Corollary~\ref{applicationcocube}
are of course not new. But the author does not know where the proofs
of these two facts can be found. So a short argumention involving
Str{\o}m's model structure is presented.

Let $\lambda$ be an ordinal. Any ordinal can be viewed as a small
category whose objects are the elements of $\lambda$, that is the
ordinal $\gamma<\lambda$, and where there exists a morphism
$\gamma\longrightarrow \gamma'$ if and only if $\gamma\leq
\gamma'$. The notation $\lambda^{op}$ will then denote the
opposite category. Let us then denote by $\C^{\lambda^{op}}$ the
category of functors from $\lambda^{op}$ to $\C$ where $\C$ is a
category. An object of $\C^{\lambda^{op}}$ is called a
\textit{tower}.

\begin{propdef}\cite{ref_model2}\cite{MR99h:55031}
Let $\C$ and $\D$ be two model categories.
A Quillen adjunction is a pair of adjoint functors
$F:\C\rightleftarrows \D:G$ between the model categories $\C$ and $\D$
such that one of the following equivalent properties holds:
\begin{enumerate}
\item if $f$ is a cofibration (resp. a trivial cofibration), then so
does $F(f)$
\item if $g$ is a fibration (resp. a trivial fibration), then so
does $G(g)$.
\end{enumerate}
One says that $F$ is a {\rm left Quillen functor}.
One says that $G$ is a {\rm right Quillen functor}.
\end{propdef}

\bd\cite{ref_model2}\cite{MR99h:55031}
An object $X$ of a model category $\C$ is {\rm cofibrant}
(resp. {\rm fibrant}) if and only if the canonical morphism
$\varnothing\longrightarrow X$ from the initial object of $\C$ to $X$
(resp. the canonical morphism $X\longrightarrow \mathbf{1}$ from $X$
to the final object $\mathbf{1}$) is a cofibration (resp.  a
fibration). \ed

\bp Let $\C$ be a model category.
There exists a model structure on $\C^{\lambda^{op}}$ such
that the limit functor $\limproj:\C^{\lambda^{op}}\longrightarrow
\C$ is a right Quillen functor and such that the fibrant towers
$T$ are exactly the towers $T:\lambda^{op}\longrightarrow \C$
such that $T_0$ is fibrant and such that for any ordinal $\gamma$
with $0\leq \gamma<\lambda$, the canonical morphism
$T_\gamma\longrightarrow \limproj_{\beta<\gamma}T_\beta$ is a
fibration of $\C$. \ep

This proposition is proved  for $\lambda=\aleph_0$ in
\cite{MR2001d:55012}.

\bpf[Sketch of proof] For a reminder about the Reedy model
structure, see \cite{ref_model2} and \cite{MR99h:55031}. With the
Reedy structure corresponding to the indexing, let us calculate
the latching space functors $L_\gamma T$ and the matching space
functors $M_\gamma T$ of a tower $T$:
\begin{enumerate}
\item if $\gamma+1<\lambda$, then $L_\gamma T=T_{\gamma+1}$
\item if $\gamma+1=\lambda$, then $L_\gamma T=\varnothing$ (the
initial object of $\C$)
\item for any $\gamma<\lambda$, $M_\gamma
T=\limproj_{\beta<\gamma}T_\beta$
\end{enumerate}
So a morphism of towers $T\longrightarrow T'$ is a cofibration
for the Reedy model structure if and only if
\begin{enumerate}
\item for $\gamma+1<\lambda$, the morphism
$T_\gamma\sqcup_{T_{\gamma+1}} T'_{\gamma+1}\longrightarrow
T'_\gamma$ is a cofibration of $\C$
\item for $\gamma+1=\lambda$, $T_\gamma\longrightarrow T'_\gamma$ is a cofibration of
$\C$.
\end{enumerate}
The limit functor $\limproj:\C^{\lambda^{op}}\longrightarrow \C$
is a right Quillen functor if and only if its left adjoint, the
constant diagram functor $\Delta:\C\longrightarrow\C^{\lambda^{op}}$ is a
left Quillen functor. Consider a cofibration $X\longrightarrow Y$
of $\C$. Then the morphism of towers $\Delta(X)\longrightarrow
\Delta(Y)$ is a cofibration if and only if either
$\gamma+1<\lambda$ and $\Delta(Y)_{\gamma+1}\longrightarrow
\Delta(Y)_\gamma$ is a cofibration or $\gamma+1=\lambda$ and
$\Delta(X)_\gamma\longrightarrow \Delta(Y)_\gamma$ is a
cofibration. This holds indeed. Therefore the limit functor is a
right Quillen functor.

But a morphism of towers $T\longrightarrow T'$ is a fibration for the
Reedy model structure if and only if for any $\gamma<\lambda$,
$T_\gamma\longrightarrow
T'_\gamma\p_{(\limproj_{\beta<\gamma}T'_\beta)}(\limproj_{\beta<\gamma}T_\beta)$
is a fibration. Hence the result. \epf

\begin{cor}\label{tourtransfinie} Let $T$ and $T'$ be two objects of
$\top^{\lambda^{op}}$ such that:
\begin{enumerate}
\item for any $\gamma<\lambda$ such that $\gamma+1<\lambda$,
the morphism $T_{\gamma+1}\longrightarrow T_\gamma$ is a
Hurewicz fibration of topological spaces
\item for any $\gamma<\lambda$ such that $\gamma$ is a limit
ordinal, the canonical morphism $T_\gamma\longrightarrow
\limproj_{\beta<\gamma}T_\beta$ is an homeomorphism
\end{enumerate}
If $f:T\longrightarrow T'$ is an objectwise homotopy
equivalence, then $\limproj f:\limproj T\longrightarrow \limproj
T'$ is a homotopy equivalence. \end{cor}

\bpf There exists a model structure on the category of topological spaces
$\top$ where the cofibration are the Hurewicz cofibrations, the
fibrations the Hurewicz fibrations and the weak homotopy equivalences
the homotopy equivalences \cite{MR35:2284} \cite{MR39:4846}
\cite{ruse}. All topological spaces are fibrant and cofibrant for this
model structure. The corollary is then due to the fact that a right
Quillen functor preserves weak homotopy equivalences between fibrant
objects and to the fact that any topological space is fibrant for this
model structure. \epf

\begin{lem} \cite{ref_model2}\cite{MR99h:55031} \label{cube} (Cube lemma)
Let $\C$ be model category. Let
\[
\xymatrix{ A_i\fr{} \fd{} & B_i \\ C_i}
\]
be two diagrams $D_i$ with $i=1,2$ of cofibrant objects of $\C$ such
that both morphisms $A_i\longrightarrow B_i$ with $i=1,2$ are
cofibrations of the model structure. Then any morphism of diagrams
$D_1\longrightarrow D_2$ which is an objectwise weak equivalence
induces a weak equivalence $\liminj D_1\longrightarrow \liminj
D_2$. \end{lem}

The dual version states as follows:

\begin{lem}\label{cocube}
Let $\C$ be a model category. Let
\[
\xymatrix{  & B_i \fd{}\\ A_i \fr{} & C_i}
\]
be two diagrams $D_i$ with $i=1,2$ of fibrant objects of $\C$ such
that both morphisms $B_i\longrightarrow C_i$ with $i=1,2$ are
fibrations of the model structure. Then any morphism of diagrams
$D_1\longrightarrow D_2$ which is an objectwise weak equivalence
induces a weak equivalence $\limproj D_1\longrightarrow \limproj
D_2$. \end{lem}

\begin{cor}\label{applicationcocube} Let
\[
\xymatrix{  & B_i \fd{}\\ A_i \fr{} & C_i}
\]
be two diagrams $D_i$ with $i=1,2$ of topological spaces
such that both morphisms $B_i\longrightarrow C_i$ with $i=1,2$
are Hurewicz fibrations. Then any morphism of
diagrams $D_1\longrightarrow D_2$ which is an objectwise homotopy
equivalence induces a homotopy equivalence
$\limproj D_1\longrightarrow \limproj D_2$. \end{cor}

\subsection{The end of the proof}

\bth \label{homotopy-chemin}
Let $X$ and $U$ be two globular complexes.  The set map
\[cat:\tgltop(X,U)\longrightarrow \tdtop(cat(X),cat(U))\] is continuous and moreover is a homotopy equivalence. \eth

\bpf  The globular decomposition of $X$ enables to view the
canonical continuous map $\varnothing\longrightarrow X$ as a
transfinite composition of $X_\beta\longrightarrow X_{\beta+1}$
for $\beta<\lambda$ such that for any ordinal $\beta<\lambda$,
one has the pushout of topological spaces
\[
\xymatrix{ \glob^{top}(\de Z_\beta)\fr{}\fd{} & X_\beta\fd{}\\
\glob^{top}( Z_\beta)\fr{} & \cocartesien X_{\beta+1}}
\]
where the pair $( Z_\beta,\de Z_\beta)$ is a NDR pair. And by
construction of the functor $cat:\gltop\longrightarrow \dtop$,
one also has for any ordinal $\beta<\lambda$ the pushout of flows
\[
\xymatrix{ \glob(\de Z_\beta)\fr{}\fd{} & cat(X_\beta)\fd{}\\
\glob( Z_\beta)\fr{} & \cocartesien cat(X_{\beta+1})}
\]
By Theorem~\ref{limittop}, one obtains the pullback of topological
spaces
\[
\xymatrix{ \tgltop(X_{\beta+1},U)\fr{}\cartesien\fd{} &
\tgltop(\glob^{top}( Z_\beta),U)\fd{} \\
\tgltop(X_{\beta},U)\fr{}& \tgltop(\glob^{top}(\de Z_\beta),U)}
\]
By Theorem~\ref{commute}, one obtains the pullback of topological
spaces
\[
\xymatrix{ \tdtop(cat(X_{\beta+1}),cat(U))\fr{}\cartesien\fd{} &
\tdtop(\glob( Z_\beta),cat(U))\fd{} \\
\tdtop(cat(X_{\beta}),cat(U))\fr{}& \tdtop(\glob(\de
 Z_\beta),cat(U))}
\]
For a given $\beta$, let us suppose that the space $\de
 Z_\beta$ is empty. Then the topological spaces $\tdtop(\glob(\de
 Z_\beta),cat(U))$ and $\tgltop(\glob^{top}(\de Z_\beta),U)$ are
both discrete. So both continuous maps
\[\tgltop(\glob^{top}( Z_\beta),U)\longrightarrow
\tgltop(\glob^{top}(\de Z_\beta),U)\] and
\[\tdtop(\glob( Z_\beta),cat(U))\longrightarrow \tdtop(\glob(\de
 Z_\beta),cat(U))\] are Hurewicz fibrations. Otherwise, if the space
 $\de Z_\beta$ is not empty, then the pair $( Z_\beta,\de Z_\beta)$ is
 a NDR pair. Then by Corollary~\ref{montage2}, the continuous map
\[\tgltop(\glob^{top}( Z_\beta),U)\longrightarrow
\tgltop(\glob^{top}(\de Z_\beta),U)\] is a Hurewicz fibration. And by
Proposition~\ref{montage}, the continuous map
\[\tdtop(\glob( Z_\beta),cat(U))\longrightarrow \tdtop(\glob(\de
 Z_\beta),cat(U))\] is a Hurewicz fibration as well. One obtains for a
 given ordinal $\beta<\lambda$ the following commutative diagram of
 topological spaces:
{\small\[
\xymatrix{ \tgltop(X_{\beta+1},U)\cartesien\fd{}\ar@{->}[rdd]
\fr{} &
\tgltop(\glob^{top}( Z_\beta),U) \ar@{->>}[d]\ar@{->}[rdd]& \\
\tgltop(X_\beta,U)\ar@{->}[rdd] \fr{} & \tgltop(\glob^{top}(\de
 Z_\beta),U) \ar@{-->}[rdd]&
\\
& \tdtop(cat(X_{\beta+1}),cat(U))\cartesien\fd{} \fr{} &
\tdtop(\glob( Z_\beta),cat(U))\ar@{->>}[d]\\
& \tdtop(cat(X_\beta),cat(U))\fr{} & \tdtop(\glob(\de
 Z_\beta),cat(U)) }
\]}
where the symbol $\xymatrix@1{\ar@{->>}[r]&}$ means Hurewicz
fibration.  One can now apply Corollary~\ref{applicationcocube}.
Therefore, if $\tgltop(X_\beta,U)\longrightarrow
\tdtop(cat(X_\beta),cat(U))$ is a homotopy equivalence of
topological spaces, then the same holds by replacing $\beta$ by
$\beta+1$. By transfinite induction, we want to prove that for
any ordinal $\beta<\lambda$, one has the homotopy equivalence of
topological spaces $\tgltop(X_\beta,U)\longrightarrow
\tdtop(cat(X_\beta),cat(U))$. The initialization is trivial: if
$X^0$ is a discrete globular complex, then $cat(X^0)=X^0$. The
passage from $\beta$ to $\beta+1$ is ensured by the proof above.
It remains to treat the case where $\beta$ is a limit ordinal.
Since the pullback of a Hurewicz fibration is a Hurewicz
fibration, then all continuous maps
\[\tgltop(X_{\beta+1},U)\longrightarrow \tgltop(X_\beta,U)\]
and
\[\tdtop(cat(X_{\beta+1}),cat(U))\longrightarrow
\tdtop(cat(X_\beta),cat(U))\] are actually Hurewicz fibrations. By Theorem~\ref{limittop}, for any limit ordinal $\beta$, one has
\[\tgltop(\liminj_{\alpha<\beta} X_\beta,U)\iso \limproj_{\alpha<\beta}\tgltop(
X_\beta,U).\] By Theorem~\ref{commute}, for any limit ordinal $\beta$,
one has
\[\tdtop(\liminj_{\alpha<\beta} cat(X_\beta),cat(U))\iso
\limproj_{\alpha<\beta}\tdtop( cat(X_\beta),cat(U)).\] The proof
is then complete with Corollary~\ref{tourtransfinie}. \epf

The preceding result can be slightly improved. The homotopy
equivalence above is actually a Hurewicz fibration. Three preliminary
propositions are necessary to establish this fact.

\bp Let $Z$ be a compact space. Let $U$ be a globular complex.
Then the canonical continuous map
\[cat:\tgltop(\glob^{top}(Z),U)\longrightarrow \tdtop(\glob(Z),cat(U))\]
is a Hurewicz fibration. \ep

\bpf 
Let $M$ be a topological space. Consider the following commutative
diagram:
\[
\xymatrix{{M \p \{0\}} \ar@{->}[r]^-{f} \ar@{^{(}->}[d]^{i}&
{\tgltop(\glob^{top}(Z),U)} \fd{cat} \\
{M \p [0,1]} \ar@{->}[r]_-{g}  \ar@{-->}[ru]^{h} &
{\tdtop(\glob(Z),cat(U))} }
\]
One has to find $h$ making the two triangles commutative where
$i:M \p \{0\}\subset M \p [0,1]$ is the canonical inclusion. Let
$h(m,u)\in\tgltop(\glob^{top}(Z),U)$ of the form
\[h(m,u)(z,t)= r(g(m,u))(z,\phi(m,z)(t))\]
where $\phi$ is a continuous map from $M\p Z$ to
$\tgltop(\vec{I}^{top},\vec{I}^{top})$. Then $cat\circ h =g$ for
any map $\phi$. It then suffices to take $\phi$ such that
\[f(m,0)(z,t)=r(g(m,u))(z,\phi(m,z)(t)).\] Such a map $\phi$ is unique
by the second assertion of Proposition~\ref{span}.  The
continuity of $\phi$ comes from its uniqueness and from the
continuity of the other components, similarly to
Proposition~\ref{span}. \epf

\bp Let $Z$ be a compact space. Let $U$ be a globular complex.
Then one has the homeomorphism
\[\tgltop(\glob^{top}(Z),U)\iso\tgltop(\glob^{top}(Z),\vI^{top})\p \tdtop(\glob(Z),cat(U)).\]
\ep

\bpf Let $f\in \tgltop(\glob^{top}(Z),U)$. Then there exists  a unique continuous
map $\phi_f:|\glob^{top}(Z)|\longrightarrow [0,1]$ such that for any
$(z,t)\in \glob^{top}(Z)$, $f(z,t)=i_U(cat(f)(z))(\phi_f(z,t))$ by
Proposition~\ref{span}. The continuous map $\phi_f$ is actually a
morphism of globular complexes from $\glob^{top}(Z)$ to
$\vI^{top}$. The mapping $f\mapsto (\phi_f,cat(f))$ defines a set map
from $\tgltop(\glob^{top}(Z),U)$ to
\[\tgltop(\glob^{top}(Z),\vI^{top})\p \tdtop(\glob(Z),cat(U))\] which is
obviously an isomorphism of sets. One obtains the isomorphism of sets
\[\tgltop(\glob^{top}(Z),U)\iso\tgltop(\glob^{top}(Z),\vI^{top})\p \tdtop(\glob(Z),cat(U)).\]
The set map
\[\tgltop(\glob^{top}(Z),\vI^{top})\p
\tdtop(\glob(Z),cat(U))\longrightarrow \tgltop(\glob^{top}(Z),U)\] is clearly
continuous for the Kelleyfication of the compact-open topology. It
remains to prove that the mapping $f\mapsto (\phi_f,cat(f))$ is
continuous. It suffices to prove that the mapping $f\mapsto \phi_f$ is
continuous since we already know that $cat(-)$ is continuous. The
latter fact comes from the continuity of the mapping $(f,z,t)\mapsto
f(z,t)$ which implies the continuity of $(f,z,t)\mapsto
\phi_f(z,t)$.  \epf

\bp Let $Z$ be a compact space. Let $U$ be a globular complex.
Let $(Z,\de Z)$ be a NDR pair.  Then the canonical continuous map
\beas
&&\tgltop(\glob^{top}(Z),U)\longrightarrow\\
&&\tgltop(\glob^{top}(\de Z),U)\p_{\tdtop(\glob(\de Z),cat(U))}
\tdtop(\glob(Z),cat(U))
\eeas
is a Hurewicz fibration. \ep

\bpf One has
\[\tgltop(\glob^{top}(Z),U)\iso \tgltop(\glob^{top}(Z),\vI^{top})
\p \tdtop(\glob(Z),cat(U))\]
and
\beas
&& \tgltop(\glob^{top}(\de Z),U)\p_{\tdtop(\glob(\de Z),cat(U))}
\tdtop(\glob(Z),cat(U))\\
&&\iso \left(\tgltop(\glob^{top}(\de Z),\vI^{top})\p \tdtop(\glob(\de Z),cat(U))\right)\\
&&\p_{\tdtop(\glob(\de Z),cat(U))}
\tdtop(\glob(Z),cat(U))\\
&&\iso\tgltop(\glob^{top}(\de Z),\vI^{top})\p \tdtop(\glob(Z),cat(U))
\eeas
So the continuous map we are studying is the cartesian product of
the Hurewicz fibration
\[\tgltop(\glob^{top}(Z),\vI^{top})\longrightarrow \tgltop(\glob^{top}(\de Z),\vI^{top})\]
by the identity of $\tdtop(\glob(Z),cat(U))$. So it is a Hurewicz fibration
as well.
\epf

\bth Let $X$ and $U$ be two globular
complexes.  The set map
\[cat:\tgltop(X,U)\longrightarrow \tdtop(cat(X),cat(U))\]
is a Hurewicz fibration. \eth

\bpf[Sketch of proof]
We use the notations of the proof of Theorem~\ref{homotopy-chemin}. We
are going to prove by transfinite induction on $\beta$ that the
canonical continuous map
\[\tgltop(X_\beta,U)\longrightarrow \tdtop(cat(X_\beta),cat(U))\]
is a Hurewicz fibration. For $\beta=0$, $X_\beta$ is the discrete
globular complex $(X^0,X^0)$. Therefore
$\tgltop(X_0,U)=\tdtop(cat(X_0),cat(U))=U^0$. Let us suppose the fact
proved for $\beta\geq 0$. Then one has the following diagram of
topological spaces
{\small\[
\xymatrix{  &
\tgltop(\glob^{top}( Z_\beta),U) \ar@{->>}[d]\ar@{->>}[rdd]& \\
\tgltop(X_\beta,U)\ar@{->>}[rdd] \fr{} & \tgltop(\glob^{top}(\de
 Z_\beta),U) \ar@{-->>}[rdd]&
\\
&  &
\tdtop(\glob( Z_\beta),cat(U))\ar@{->>}[d]\\
& \tdtop(cat(X_\beta),cat(U))\fr{} & \tdtop(\glob(\de
 Z_\beta),cat(U)) }
\]}
where the symbol $\xymatrix@1{\ar@{->>}[r]&}$ means Hurewicz
fibration.  We then consider the Reedy category
\[\xymatrix{
& 2\fd{} \\
0\fr{}& 1}
\]
and the Reedy model category of diagrams of topological spaces over
this small category. In this model category, the fibrant diagrams $D$
are the diagrams such that $D_0$, $D_1$ and $D_2$ are fibrant and such
that $D_2\longrightarrow D_1$ is a fibration. And a morphism of
diagrams $D\longrightarrow D'$ is fibrant if and only if both
$D_0\longrightarrow D'_0$ and $D_1\longrightarrow D'_1$ are fibrant
and if $D_2\longrightarrow D_1\p_{D'_1}D'_2$ is fibrant. So it remains
to check that the inverse limit functor is a right Quillen functor to
complete the proof. It then suffices to prove that the constant
diagram functor is a left Quillen functor. For this Reedy model
structure, a morphism of diagrams $D\longrightarrow D'$ is cofibrant
if both $D_0\longrightarrow D'_0$ and $D_2\longrightarrow D'_2$ are
cofibrant and if $D'_0\sqcup_{D_0}D_1\longrightarrow D'_1$ is
cofibrant. So a diagram $D$ is cofibrant if and only if $D_0$, $D_1$
and $D_2$ are cofibrant and if $D_0\longrightarrow D_1$ is a
cofibration. Hence the result.  \epf

\begin{cor} Let $X$ and $U$ be two globular
complexes.  The set map
\[cat:\tgltop(X,U)\longrightarrow \tdtop(cat(X),cat(U))\]
is onto. \end{cor}

\bpf Any Hurewicz fibration which is a homotopy equivalence is onto since it satisfies
the right lifting property with respect to $\varnothing\longrightarrow \{0\}$. \epf

\section{Comparison of S-homotopy in $\gltop$ and in
$\dtop$}\label{compS}

\subsection{Pairing $\boxtimes$ between a topological space and a
flow}

\bd\cite{model3} Let $U$ be a topological space. Let $X$ be a
flow. The flow $\{U,X\}_S$ is defined as follows:
\begin{enumerate}
\item The $0$-skeleton of $\{U,X\}_S$ is $X^0$.
\item For $\alpha,\beta\in X^0$, the topological space ${\P}_{\alpha,\beta}\{U,X\}_S$
is $\ttop(U,{\P}_{\alpha,\beta}X)$ with an obvious definition of
the composition law.
\end{enumerate}
\ed

\bth \label{adj} (\cite{model3} Theorem~7.8) \label{existenceflow} 
Let $U$ be a topological space. The functor $\{U,-\}_S$ has a left adjoint
which will be denoted by $U\boxtimes -$. Moreover: 
\begin{enumerate}
\item one has the natural isomorphism of flows \[U\boxtimes (\liminj_i X_i)
\iso \liminj_i (U\boxtimes X_i)\]
\item there is a natural isomorphism of flows $\{*\}\boxtimes Y\iso Y$
\item if $Z$ is another topological space, one has the natural
isomorphism of flows \[U\boxtimes \glob(Z)\iso
\glob(U\p Z)\]
\item for any flow $X$ and any topological space $U$, one has
the natural bijection of sets \[(U\boxtimes X)^0\iso
X^0\]
\item if $U$ and $V$ are two topological spaces, then $(U\p V)\boxtimes Y\iso
U\boxtimes (V\boxtimes Y)$ as flows
\item for any flow $X$, $\varnothing \boxtimes X\iso X^0$.
\end{enumerate}
\eth

\subsection{S-homotopy of flows}

\bd \cite{model3} A morphism of flows $f:X\longrightarrow Y$ is
said \textit{synchronized} if and
only if it induces a bijection of sets between the
$0$-skeleton of $X$ and the
$0$-skeleton of $Y$. \ed

\bd \cite{model3} Two morphisms of flows $f$ and $g$ from $X$ to
$Y$ are {\rm S-homotopy equivalent} if and only if there exists
\[H\in\top([0,1],\tdtop(X,Y))\] such that $H(0)=f$ and $H(1)=g$. We
denote this situation by $f\sim_S g$. \ed

\bd \cite{model3} Two flows are {\rm S-homotopy equivalent} or 
{\rm S-homotop\-ic} if and only if there exist morphisms of flows
$f:X\longrightarrow Y$ and $g:Y\longrightarrow X$ such that $f\circ
g\sim_S \id_Y$ and $g\circ f\sim_S \id_X$. \ed

\bp (\label{caracflow} Proposition~7.5) 
\cite{model3} Let $f$ and $g$ be two
morphisms of flows from $X$ to $Y$. Then $f$ and $g$ are S-homotopy
equivalent if and only if there exists a continuous map \[h\in
\top([0,1],\tdtop(X,Y))\] such that
$h(0)=f$ and $h(1)=g$. \ep

\bp (\cite{model3} Corollary~7.11) [Cylinder functor] \label{Scyl} The mapping
$X\mapsto [0,1]\boxtimes X$ induces a functor from $\dtop$ to
itself which is a cylinder functor with
the natural transformations $e_i:\{i\}\boxtimes - \to
[0,1]\boxtimes -$ induced by the inclusion maps $\{i\}\subset
[0,1]$ for $i\in\{0,1\}$ and with the natural transformation
$p:[0,1]\boxtimes -\longrightarrow \{0\}\boxtimes -$ induced by the
constant map $[0,1]\longrightarrow \{0\}$. Moreover, two morphisms of
flows $f$ and $g$ from $X$ to $Y$ are S-homotopic if and only if
there exists a morphism of flows $H:[0,1]\boxtimes X\to
Y$ such that $H\circ e_0=f$ and $H\circ e_1=g$. Moreover
$e_0\circ H\sim_{S} Id$ and $e_1\circ H\sim_{S} Id$. \ep

\subsection{Pairing $\boxtimes$ and S-homotopy}

\bp\label{reduc1} \label{cylindre-path}
Let $U$ be a compact space. Let $X$ be a globular complex. Then one has
the isomorphism of flows $cat(U\boxtimes X)\iso U\boxtimes
cat(X)$. \ep

\bpf Let $(\de Z_\beta, Z_\beta,\phi_\beta)_{\beta<\lambda}$ be the globular
decomposition of $X$.  This is clear if $X=X_0=(X^0,X^0)$ and if $X=\glob^{top}(Z)$ where
$Z$ is compact. It then suffices to make a transfinite induction on $\beta$
to prove $cat(U\boxtimes X_\beta)\iso U\boxtimes cat(X_\beta)$. \epf

\bth \label{cor} \label{corS} The set map
$cat:\gltop(X,U)\to
\dtop(cat(X),cat(U))$ induces a
bijection of sets $\gltop(X,U)/\!\!\sim_{S}\iso
\dtop(cat(X),cat(U))/\!\!\sim_{S}$.
\eth

\bpf Let $f$ and $g$ be two S-homotopy equivalent morphisms of
globular complexes from $X$ to $Y$. Then there exists a morphism of
globular complexes $H:[0,1]\boxtimes X\longrightarrow Y$ such that the
composite $H\circ e_0$ is equal to $f$ and the composite $H\circ e_1$
is equal to $g$. Then $cat(H):[0,1]\boxtimes X\longrightarrow Y$
induces by Proposition~\ref{cylindre-path} a S-homotopy between
$cat(f)$ and $cat(g)$. So the mapping $cat$ induces a set map
$\gltop(X,U)/\!\!\sim_{S}\to \dtop(cat(X),cat(U))/\!\!\sim_{S}$. By
Proposition~\ref{caracgl}, the set $\gltop(X,U)/\!\!\sim_{S}$ is
exactly the set of path-connected components of $\tgltop(X,U)$. By
Proposition~\ref{caracflow}, the set
$\dtop(cat(X),cat(U))/\!\!\sim_{S}$ is exactly the set of
path-connected components of $\tdtop(cat(X),cat(U))$. But the set map
$cat:\tgltop(X,U)\to \tdtop(cat(X),cat(U))$ induces a homotopy
equivalence by Theorem~\ref{homotopy-chemin}. So the two topological
spaces have the same path-connected components. \epf

\begin{cor} Two globular complexes are
S-homotopy equivalent if and only if the
corresponding flows are S-homotopy equivalent.
\end{cor}

\begin{cor} The localization of the category of globular complexes
with respect to the class of S-homotopy equivalences is equivalent to
the localization of the full and faithful subcategory of flows of the
form $cat(X)$ with respect to the S-homotopy equivalences. \end{cor}

\bpf 
This is due to the existence of the cylinder functor both for the
S-homotopy of globular complexes and for the S-homotopy of flows. \epf

\section{Conclusion}

This part shows that the category of flows is an appropriate framework
for the study of S-homotopy equivalences. The category $\dtop$ has
nicer categorical properties than $\gltop$, for example because it is
both complete and cocomplete.

\part{Flow up to weak S-homotopy}\label{model4}

\section{Introduction}

We prove that the functor from the category of globular CW-complexes
to the category of flows induces an equivalence of categories from the
localization of the category of globular CW-complexes with respect to
the class of the S-homotopy equivalences to the localization of the
category of flows with respect to the class of weak S-homotopy
equivalences.

\section{The model structure of $\dtop$}

Some useful references for the notion of model category are
\cite{MR99h:55031} \cite{MR2001d:55012}. See also \cite{ref_model1} 
\cite{ref_model2}.

\bth (\cite{model3} Theorem~19.7) \label{wW}  The category of flows
can be given a model structure such that:
\begin{enumerate}
\item The weak equivalences  are the {\rm weak S-homotopy equivalences}, that is
a morphism of flows $f:X\longrightarrow Y$ such that
$f:X^0\longrightarrow Y^0$ is an isomorphism of sets and $f:{\P}
X\longrightarrow {\P} Y$ a weak homotopy equivalence of topological
spaces.
\item The fibrations are the continuous maps satisfying the RLP
 with respect to the morphisms $\glob(\mathbf{D}^n)\longrightarrow
 \glob([0,1]\p \mathbf{D}^n)$ for $n\geq 0$. The fibrations are
 exactly the morphisms of flows $f:X\longrightarrow Y$ such that ${\P}
 f:{\P} X\longrightarrow {\P} Y$ is a Serre fibration of $\top$.
\item The cofibrations are the morphisms satisfying the LLP  with
respect to any map satisfying the RLP with respect to the morphisms
$\glob(\mathbf{S}^{n-1})\longrightarrow \glob(\mathbf{D}^n)$ for
$n\geq 0$ and with respect to the morphisms
$\varnothing\longrightarrow
\{0\}$ and $\{0,1\}\longrightarrow \{0\}$.
\item Any flow is fibrant.
\end{enumerate}
\eth

\begin{nota} Let $\mathcal{S}$ be the subcategory of weak S-homotopy
equivalences. Let $I^{gl}$ be the set of morphisms of flows
$\glob(\mathbf{S}^{n-1})\longrightarrow \glob(\mathbf{D}^n)$ for $n\geq 0$. Let
$J^{gl}$ be the set of morphisms of flows
$\glob(\mathbf{D}^{n})\longrightarrow \glob([0,1]\p \mathbf{D}^n)$. Notice that all
arrows of  $\mathcal{S}$, $I^{gl}$ and $J^{gl}$ are
synchronized. At last, denote by $I^{gl}_+$ be the union of
$I^{gl}$ with the two morphisms of flows $R:\{0,1\}\longrightarrow
\{0\}$ and $C:\varnothing\subset \{0\}$. \end{nota}

\section{Strongly cofibrant replacement of a flow}

\bd Let $X$ be a flow. Let $n\geq 0$. Let
$f_i:\glob(\mathbf{S}^{n-1})\longrightarrow X$ be a family of
morphisms of flows with $i\in I$ where $I$ is some set.  Then the
pushout $Y$ of the diagram
\[
\xymatrix{ \bigsqcup_{i\in I}
\glob(\mathbf{S}^{n-1})\ar@{->}[rr]^-{\bigsqcup_{i\in
I} f_i } \fd{\subset} &&  X \\ \bigsqcup_{i\in I}
\glob(\mathbf{D}^n)&&}
\]
is called a {\rm $n$-globular extension} of $X$. The family of
$f_i:\glob(\mathbf{S}^{n-1})\longrightarrow X$ is called the {\rm
globular decomposition of the extension}. \ed

\bd Let $i:A\longrightarrow X$ be a morphism of flows. Then the
morphism $i$ is a {\rm relative globular extension} if the flow $X$ is
isomorphic to a flow $X_\omega=\liminj X_n$ such that for any integer
$n\geq 0$, $X_n$ is a $n$-globular extension of $X_{n-1}$ (by
convention, let $X_{-1}=A$). One says that $\dim(X,A)= p$ if
$X_\omega=X_p=X_{p+1}=\dots$ and if $X_{p-1}\neq X_p$. The flow $X_n$
is called the $n$-skeleton of $(X,A)$ and the family of $(X_n)_{n\geq
0}$ the skeleton. \ed

\bd 
A flow $X$ is said {\rm strongly cofibrant} if and only if the pair
$(X,X^0)$, where $X^0$ is the $0$-skeleton, is a relative globular
extension. Let \[\dim(X)=\dim(X,X^0).\] \ed

Notice that any strongly cofibrant flow is cofibrant for the model
structure of $\dtop$. Using Theorem~\ref{wW}, we already know that any
flow is weakly S-homotopy equivalent to a cofibrant flow and that this
cofibrant flow is unique up to S-homotopy.  Such a cofibrant flow is
usually called a \textit{cofibrant replacement}.  With the standard
construction of the cofibrant replacement involving the ``Small Object
Argument'', we can only say that the cofibrant replacement of a flow
can be taken in the $I^{gl}_+$-cell complexes.

We want to prove in this section that the cofibrant replacement can be
supposed strongly cofibrant. This is therefore a stronger statement
than the usual one.

\bth\label{classiquen}\label{comp2} (\cite{model3} Theorem~15.2) 
Suppose that one has the pushout of flows
\[
\xymatrix{
\glob(\mathbf{S}^{n}) \fr{}\fd{} & A \fd{} \\
\glob(\mathbf{D}^{n+1}) \fr{} & \cocartesien X }
\]
for some $n\geq 1$. Then the continuous map $f:{\P} A\longrightarrow
{\P} X$ is a closed $n$-connected inclusion. \eth

\bth \label{image} 
Any flow is weakly S-homotopy equivalent to a strongly cofibrant
flow. This "strongly cofibrant replacement" is unique up to
S-homotopy. \eth

\bpf As usual in this kind of proof, two kinds of processes are
involved; the first is that of attaching cells like
$\glob(\mathbf{S}^n)$ so as to create new generators; the second,
of attaching cells like $\glob(\mathbf{D}^n)$ to create new
relations.

Let $X$ be an object of $\dtop$. Let $T_{-1}=X^0$ (so ${\P}
T_{-1}=\varnothing$). Then the canonical morphism
$f_{-1}:T_{-1}\longrightarrow X$ is synchronized.  If ${\P}
X=\varnothing$, then the proof is ended. Otherwise, for any $\gamma\in
{\P} X$, let us attach a copy of $\vI$ such that $[0,1]\in {\P} \vI$
is mapped to $\gamma$. Then the canonical morphism of flows
$f_0:T_0\longrightarrow X$ induces an onto map $\pi_0(f_0):\pi_0({\P}
T_0)\longrightarrow \pi_0({\P} X)$ (where $\pi_i(U)$ is the $i$-th
homotopy group of $U$). In other terms, $T_0$ is the flow having $X^0$
as $0$-skeleton and the set ${\P} X$ equipped with the discrete
topology as path space.

We are going to introduce by induction on $n\geq 0$ a $n$-globular
extension $T_n$ of $T_{n-1}$ such that the canonical morphism of flows
$f_n:T_n\longrightarrow X$ satisfies the following conditions:
\begin{enumerate}
\item the morphism of flows $f_n$ is synchronized
\item for any base-point $\gamma$, $\pi_n(f_n):\pi_n({\P} T_n,\gamma)\longrightarrow \pi_n({\P} X,\gamma)$
is onto
\item for any base-point $\gamma$, and for any $0\leq i<n$,
$\pi_i(f_n):\pi_i({\P} T_n,\gamma)\longrightarrow \pi_i({\P}
X,\gamma)$ is an isomorphism.
\end{enumerate}

The passage from $T_0$ to $T_1$ is fairly different from the rest of
the induction. To obtain a bijection $\pi_0(f_1):\pi_0({\P}
T_1)\longrightarrow \pi_0({\P} X)$, it suffices to have a bijection
$\pi_0(f_1):\pi_0({\P}_{\alpha,\beta} T_1)\longrightarrow
\pi_0({\P}_{\alpha,\beta} X)$ for any $\alpha,\beta\in X^0$. Let $x$
and $y$ be two distinct elements of $\pi_0({\P}_{\alpha,\beta} T_0)$
having the same image in $\pi_0({\P}_{\alpha,\beta} X)$. Then $x$ and
$y$ correspond to two non-constant execution paths $\gamma_x$ and
$\gamma_y$ from $\alpha$ to $\beta$. Consider the morphism of flows
$\glob(\mathbf{S}^0)\longrightarrow X$ such that $-1\mapsto
\gamma_x$ and $1\mapsto \gamma_y$. Then let us attach a cell
$\glob(\mathbf{D}^1)$ by the pushout
\[
\xymatrix{
\glob(\mathbf{S}^0)\fr{} \fd{} & T_0 \fd{}\ar@/^15pt/[rdd]&\\
\glob(\mathbf{D}^1) \ar@/_15pt/[rrd]_{g}\fr{} & \cocartesien T_0^{(1)}\ar@{-->}[rd]^{k^{(1)}} &\\
&&X}
\]
By construction, the equality $x=y$ holds in $T_0^{(1)}$. By
transfinite induction, one obtains a flow $U_0$ and a morphism of
flows $U_0\longrightarrow X$ inducing a bijection $\pi_0(U_0)\iso
\pi_0(X)$. We now have to make $\pi_1(U_0)\longrightarrow
\pi_1(X)$ onto. The passage from $U_0$ to $T_1$ is analogous to
the passage from $U_n$ to $T_{n+1}$ for $n\geq 1$, as explained below.

Let us suppose $T_n$ constructed for $n\geq 1$. We are going to
construct the morphism $T_n\longrightarrow T_{n+1}$ as a transfinite
composition of pushouts of the morphism of flows
$\glob(\mathbf{S}^n)\longrightarrow \glob(\mathbf{D}^{n+1})$.  By
Theorem~\ref{comp2}, the pair $({\P} T_{n+1},{\P} T_n)$ will be
$n$-connected, and so the canonical maps $\pi_i({\P}
T_n)\longrightarrow \pi_i({\P} T_{n+1})$ will be bijective for $i< n$.
So the canonical map $\pi_i({\P} T_{n+1})\longrightarrow
\pi_i({\P} X)$ will remain bijective for $i< n$. By induction
hypothesis, the map $\pi_n(f_n):\pi_n({\P} T_n,\gamma)\longrightarrow
\pi_n({\P} X,\gamma)$ is onto. To each element of $\pi_n({\P} T_n,\gamma)$
with trivial image in $\pi_n({\P} X,\gamma)$ corresponds a continuous
map $\mathbf{S}^n\longrightarrow {\P} T_n$. Since $\mathbf{S}^n$ is
connected, it can be associated to a morphism of flows
$\glob(\mathbf{S}^n)\longrightarrow T_n$. Let us attach to $T_n$ a
cell $\glob(\mathbf{D}^{n+1})$ using the latter morphism. And repeat
the process transfinitely. Then one obtains a relative
$(n+1)$-globular extension $U_n$ of $T_n$ such that $\pi_i({\P}
U_n)\longrightarrow \pi_i({\P} X)$ is still bijective for $i< n$ 
and such that $\pi_n({\P} U_n)\longrightarrow \pi_n({\P} X)$
becomes bijective. Now we have to make $\pi_{n+1}({\P}
U_n,\gamma)\longrightarrow \pi_{n+1}({\P} X,\gamma)$ onto for any
base-point $\gamma$. Let
$g:(\mathbf{D}^{n+1},\mathbf{S}^n)\longrightarrow ({\P} X,\gamma)$ be a
relative continuous map which corresponds to an element of
$\pi_{n+1}({\P} X)$. Let us consider the following commutative diagram:
\[
\xymatrix{
\glob(\mathbf{S}^n)\fr{\gamma_*} \fd{} & U_n \fd{} \ar@/^15pt/[rdd]& \\
\glob(\mathbf{D}^{n+1}) \ar@/_15pt/[rrd]_{g}\fr{} &
U_n^{(1)}\cocartesien
\ar@{-->}[rd]^{k^{(1)}} & \\
& & X }
\]
where $\gamma_*(0)=s(\gamma)$, $\gamma_*(1)=t(\gamma)$ and for any
$z\in \mathbf{S}^n$, $\gamma_*(z)=\gamma$. Then because of the
universal property satisfied by the pushout, there exists a morphism
of flows $k^{(1)}:U_n^{(1)}\longrightarrow X$ and by construction, the
canonical morphism $\mathbf{D}^{n+1}\longrightarrow {\P} U_n^{(1)}$ is
an inverse image of $g$ by the canonical map $\pi_{n+1}({\P}
U_n^{(1)},\gamma)\longrightarrow \pi_{n+1}({\P} X,\gamma)$. By
transfinite induction, one then obtains for some ordinal $\lambda$ a
flow $U_n^{(\lambda)}$ such that $\pi_{n+1}({\P}
U_n^{(\lambda)},\gamma)\longrightarrow \pi_{n+1}({\P} X,\gamma)$ is
onto.  It then suffices to set $T_{n+1}:=U_n^{(\lambda)}$. The colimit
$\liminj T_n$ is then a strongly cofibrant replacement of $X$ and
$\liminj f_n: \liminj T_n\longrightarrow X$ is then a weak S-homotopy
equivalence by construction. The uniqueness of this strongly cofibrant
replacement up to S-homotopy is a consequence of
Theorem~\ref{wW}. \epf

\section{The category of S-homotopy types}\label{consS}

\bth\label{equi} The functor $cat$ from $\gltop$ to $\dtop$
induces an equivalence between the localization
$\glCW[\mathcal{SH}^{-1}]$ of globular CW-complexes with respect to
the class $\mathcal{SH}$ of S-homotopy equivalences and the
localization of the full and faithful subcategory of $\dtop$
consisting of the strongly cofibrant flows by the S-homotopy
equivalences. \eth

\bpf Let $X$ be a strongly cofibrant flow. Let
$(X_n)_{n\geq 0}$ be the skeleton of the relative globular
extension $(X,X^0)$ (with the convention $X_{-1}=X^0$). Let
$P(n)$ be the statement: ``there exists a
globular CW-complex $Y$ of dimension $n$
such that $cat(Y)=X_n$ (by convention a globular CW-complex
of dimension $-1$ will be a discrete space)''. Suppose
$P(n)$ proved for $n\geq -1$. Using Theorem~\ref{i}, choose
a continuous map $i_{Y}:{\P} Y\longrightarrow {\P}^{top}Y$. Let
\[
\xymatrix{ \bigsqcup_{i\in I}
\glob(\mathbf{S}^{n})\ar@{->}[rr]^-{\bigsqcup_{i\in I}
f_i } \fd{\subset} &&  X_n\fd{} \\ \bigsqcup_{i\in I}
\glob(\mathbf{D}^{n+1})\ar@{->}[rr]&&X_{n+1}\cocartesien}
\]
be the pushout defining $X_{n+1}$. Then the pushout of multipointed
spaces
\[
\xymatrix{ \bigsqcup_{i\in I}
\glob^{top}(\mathbf{S}^{n})\ar@{->}[rrrr]^-{\bigsqcup_{i\in
I} (z,t) \mapsto
i_Y(f_i(z))(t) } \fd{\subset}
&&&&  Y\fd{} \\ \bigsqcup_{i\in I}
\glob^{top}(\mathbf{D}^{n+1})\ar@{->}[rrrr]&&&&Y'\cocartesien}
\]
gives the solution. It remains to prove that the functor is both full
and faithful. Since S-homotopy in $\glCW$ is characterized by a
cylinder functor (cf. \cite{diCW} or Corollary~\ref{cylgl}), one has
the natural bijection of sets
\[\glCW[\mathcal{SH}^{-1}](X,Y)\iso
\gltop(X,Y)/\!\!\sim_S\] for any globular CW-complexes
$X$ and $Y$.  Since S-homotopy in $\dtop$ is also characterized
by a cylinder functor (cf. Proposition~\ref{Scyl}), one also has the
natural bijection of sets
\[\dtop(cat(X),cat(Y))/\!\!\sim_S\iso
\dtop[\mathcal{SH}^{-1}](cat(X),cat(Y)).\] The
theorem is then a consequence of Theorem~\ref{corS}.  \epf

\bth \label{equiS} 
The localization $\dtop[\mathcal{S}^{-1}]$ of $\dtop$ with respect to
the class $\mathcal{S}$ of weak S-homotopy equivalences exists (i.e.
is locally small). The functor $cat:\glCW\longrightarrow \dtop$
induces an equivalence of categories $\glCW[\mathcal{SH}^{-1}]\iso
\dtop[\mathcal{S}^{-1}]$. \eth

\bpf Let $X$ be an object of $\dtop$. By Theorem~\ref{image},
there exists a strongly cofibrant flow $X'$ weakly
S-homotopy equivalent to $X$. By
Theorem~\ref{equi}, there exists a globular
CW-complex $Y$ with
$cat(Y)\iso X'$.  So $cat(Y)$ is isomorphic
to $X$ in $\dtop[\mathcal{S}^{-1}]$. So the functor 
$cat:\glCW[\mathcal{SH}^{-1}]\longrightarrow
\dtop[\mathcal{S}^{-1}]$ is essentially surjective.

Let $Y_1$ and $Y_2$ be two globular CW-complexes. Then
\[\glCW[\mathcal{SH}^{-1}](Y_1,Y_2)\iso
\dtop[\mathcal{SH}^{-1}](cat(Y_1),cat(Y_2))\iso
\dtop[\mathcal{S}^{-1}](cat(Y_1),cat(Y_2))\] the last isomorphism
being due to the facts that $cat(Y_1)$ is cofibrant and that
$cat(Y_2)$ is fibrant for the model structure of $\dtop$. Therefore
$cat:\glCW[\mathcal{SH}^{-1}]\longrightarrow
\dtop[\mathcal{S}^{-1}]$ is full and faithful. \epf

\begin{cor} Let $\CW$ be the category of CW-complexes. Let
$\top$ be the category of compactly generated topological spaces.  Let
$\ho(\CW)$ be the localization of $\CW$ with respect to homotopy
equivalences and $\ho(\top)$ be the localization of $\top$ with
respect to weak homotopy equivalences. Then the commutative diagram
\[
\xymatrix{ {\CW} \fr{} \fd{\glob^{top}(-)} & {\top}\fd{\glob(-)}\\
{\glCW} \fr{} &  {\dtop}}
\]
gives rise to the commutative diagram
\[
\xymatrix{ {\ho(\CW)} \fr{\iso}
\ar@{^{(}->}[d]_{\glob^{top}(-)} &
{\ho(\top)}
\ar@{^{(}->}[d]^{\glob(-)}\\
{\glCW[\mathcal{SH}^{-1}]} \ar@{->}[r]^{\iso} &
{\dtop[\mathcal{S}^{-1}]}}
\]
\end{cor}

\section{Conclusion}

The model structure of \cite{model3} on the category of flows provides
a new interpretation of the notion of S-homotopy equivalence. It
allowed us to prove in Part~\ref{model4} that the functor from the
category of globular CW-complexes to the category of flows induces an
equivalence of categories from the localization of the category of
globular CW-complexes with respect to the class of the S-homotopy
equivalences to the localization of the category of flows with respect
to the class of weak S-homotopy equivalences.

\part{T-homotopy and flow}\label{model2.5}

\section{Introduction}

The purpose of this part is the construction of a class of morphisms
of flows, the T-homotopy equivalences, so that the following theorem
holds:

\bth 
Let $X$ and $U$ be globular complexes. If $f:X\longrightarrow U$ is a
T-homotopy equivalence of globular complexes, then
$cat(f):cat(X)\longrightarrow cat(U)$ is a T-homotopy equivalence of
flows. Conversely, if $g:cat(X)\longrightarrow cat(U)$ is a T-homotopy
equivalence of flows, then $g=cat(f)$ for some T-homotopy equivalence
of globular complexes $f:X\longrightarrow U$. \eth

where:

\bd A \textit{T-homotopy} is a 
morphism $f:X\longrightarrow Y$ of globular complexes inducing an
homeomorphism between the two underlying topological spaces. \ed

Section~\ref{sectionT} defines the class of T-homotopy equivalences in
the category of flows. Section~\ref{compT} is devoted to proving the
theorem above.

\section{T-homotopy in $\dtop$}\label{sectionT}

\begin{figure}
\begin{center}
\includegraphics[width=5cm]{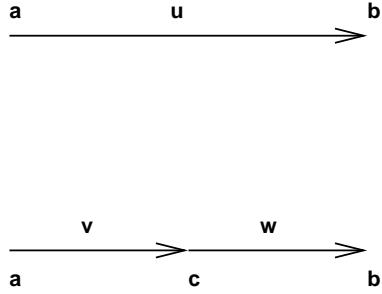}
\end{center}
\caption{Concatenation of $v$ and $w$} \label{ex1}
\end{figure}

The idea of T-homotopy is to change nothing globally except that new
states may appear in the middle of full globes. In particular, the
additional states appearing in the $0$-skeleton must not create any
new branching or merging areas of execution paths. For example, the
unique morphism of flows $F$ such that $F(u)=v* w$ in Figure~\ref{ex1}
is a T-homotopy.

\begin{figure}
\begin{center}
\includegraphics[width=5cm]{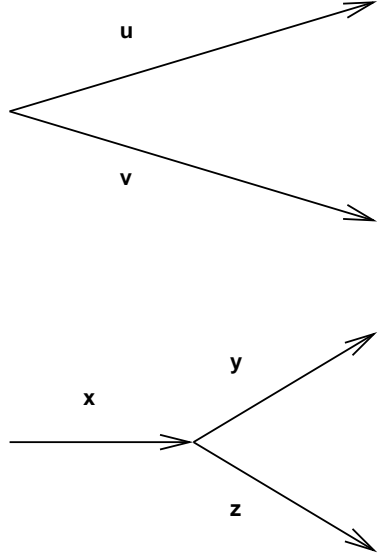}
\end{center}
\caption{Two $1$-dimensional automata not T-homotopy equivalent}
\label{ex2}
\end{figure}

We need again the notion of quasi-flow introduced in Part~\ref{model1}
Section~\ref{quasi}. Recall that a flow can be viewed as a particular
case of quasi-flow.

\bd Let $X$ be a quasi-flow. Let $Y$ be a subset of $X^0$. Then
the \textit{restriction} $X\!\restriction_Y$ of $X$ over $Y$ is
the unique quasi-flow such that
$(X\!\restriction_Y)^0=Y$ and such that
\[{\P}^{top} (X\!\restriction_Y) = \bigsqcup_{(\alpha,\beta)\in Y\p Y}
{\P}^{top}_{\alpha,\beta} X\] equipped with the topology induced
by the one of ${\P}^{top} X$. \ed

Let $X$ be a flow. As in \cite{fibrantcoin}, let $\mathcal{R}^-$ be
the smallest closed equivalence relation on ${\P} X$ identifying
$\gamma_1$ and $\gamma_1*\gamma_2$ whenever $\gamma_1$ and
$\gamma_1*\gamma_2$ are defined in ${\P} X$, and let \[{\P}^- X= {\P}
X/ \mathcal{R}^-.\] Symmetrically, let us consider the smallest closed
equivalence relation $\mathcal{R}^+$ identifying $\gamma_2$ and
$\gamma_1* \gamma_2$ if $\gamma_2$ and $\gamma_1* \gamma_2$ belong to
${\P} X$. Then let \[{\P}^+ X= {\P} X/ \mathcal{R}^+.\]

\bd\label{Tdi} A morphism of flows $f:X\longrightarrow Y$ is a T-homotopy
equivalence if and only if the following conditions are satisfied:
\begin{enumerate}
\item The morphism of flows $f:X\to
Y\!\restriction_{f(X^0)}$ is an isomorphism of flows. In
particular, the set map $f^0:X^0\longrightarrow Y^0$ is one-to-one.
\item For any $\alpha\in Y^0\backslash f(X^0)$, the topological spaces
${\P}^-_\alpha Y$ and ${\P}^+_\alpha Y$ are singletons.
\item For any $\alpha\in Y^0\backslash f(X^0)$, there are
execution paths $u$ and $v$ in $Y$ such that $s(u)\in f^0(X^0)$,
$t(u)=y$, $s(v)=y$ and $t(v)\in f^0(X^0)$.
\end{enumerate}
\ed

The first condition alone does not suffice for a characterization of
T-homotopy, since the unique morphisms of flows $F'$ such that
$F'(u)=v$ satisfies this condition as well. The additional state (i.e.  
$b=tw$) creates a new final state.

Now consider Figure~\ref{ex2}. In the globular complex setting, there
are no T-homoto\-py equivalences between them because the underlying
topological spaces are not homeomorphic because of the calculation $x$
before the branching.  However the unique morphisms of flows $F$ such
that $F(u)=x*y$ and $F(v)=x*z$ satisfies the first and third
conditions of Definition~\ref{Tdi}, but not the second one.

Requiring that ${\P}_\alpha^- Y$ and ${\P}_\alpha^+ Y$ are only
contractible for $\alpha\in Y^0 \backslash f(X^0)$ is not sufficient
either. Indeed, consider two contractible topological spaces $X$ and
$Y$ and the morphism of globular complexes $f:\glob^{top}(X\p Y)\to
\glob^{top}(X)*\glob^{top}(Y)$ such that $f((x,y),t)=(x,2t)$ for 
$0\leq t\leq 1/2$ and $f((x,y),t)=(y,2t-1)$ for $1/2\leq t\leq 1$. The
morphism of flows $f$ would be a T-homotopy equivalence.

The third condition is also necessary because otherwise, the directed
segment $\vI$ would be T-homot\-opy equivalent to the disjoint sum of
$\vI$ with the concatenation of an infinite number of copies of $\vI$.

\section{Comparison of T-homotopy in $\gltop$ and in
$\dtop$}\label{compT}

\subsection{Properties of T-homotopy}

Some useful properties of T-homotopy equivalences of flows are proved
in this section.

\bth\label{injecte} Let $f$ be a morphism of flows from $X$ to
$Y$. Assume that $f$ is the pushout of a morphism of flows of the form
$cat(g):cat(U)\longrightarrow cat(V)$ where $g:U\longrightarrow V$ is
a T-homotopy equivalence of globular complexes. Then the morphism of
flows $X\longrightarrow Y\!\restriction_{f(X^0)}$ is an isomorphism of
flows. In particular, the continuous map ${\P} f:{\P} X\longrightarrow
{\P} Y$ is one-to-one. \eth

\bpf First of all, assume that $X=cat(U)$, $Y=cat(V)$ and $f=cat(g)$ for
some T-homotopy equivalence $g:U\longrightarrow V$. The morphism of
quasi-flows $qcat(g):qcat(U)\longrightarrow
qcat(V)\!\restriction_{g(X^0)}$ has an obvious inverse from
$qcat(V)\!\restriction_{g(X^0)}$ to $qcat(U)$ denoted by
$qcat(g)^{-1}$ sending $\gamma\in
{\P}^{top}qcat(V)\!\restriction_{g(X^0)}$ to $g^{-1}\circ
\gamma\in {\P}^{top} qcat(U)$. Using the natural transformation
$p:qcat\longrightarrow cat$, one obtains that \[p\left(qcat(g)^{-1}\right):{\P}
V\!\restriction_{g(X^0)}\longrightarrow {\P} U\] is an inverse continuous map
of ${\P} g:{\P} U\longrightarrow {\P} V\!\restriction_{g(X^0)}$.

Now take a general T-homotopy equivalence of flows $f$ from $X$ to
$Y$. By hypothesis, there exists a cocartesian diagram of flows
\[
\xymatrix{
cat(U) \fr{}\fd{cat(g)} & X \fd{}\\
cat(V) \fr{} & \cocartesien Y}
\]
for some T-homotopy equivalence of globular complexes
$f:U\longrightarrow V$.  Consider the following commutative diagram of
flows
\[
\xymatrix{
cat(U)\fr{} \fd{cat(g)} & X \fd{}\ar@/^15pt/[rdd]^{\phi_1}&\\
cat(V)\!\restriction_{g(X^0)} \ar@/_15pt/[rrd]_{\phi_2}\fr{} & Y\!\restriction_{f(X^0)}\ar@{-->}[rd]^{h} &\\
&&Z}
\]
One wants to prove the existence of $h$ making the diagram
commutative. So consider the following diagram (where a new flow
$\overline{Z}$ is defined as a pushout)
\[
\xymatrix{ && cat(U) \ar@{->}[dl]_{cat(g)}\fr{}\fd{cat(g)} & X
\ar@{->} `u[l]  `[lll] `[dd]_{\phi_1} [lldd]
\fd{}\\
&cat(V)\!\restriction_{g(X^0)} \fd{\phi_2}\ar@{^{(}->}[r]& cat(V) \fd{\overline{\phi_2}} \fr{} & \cocartesien Y\\
&Z \fr{u} & \cocartesien \overline{Z}& }
\]
Every part of this diagram is commutative. Therefore one obtains the
commutative diagram
\[
\xymatrix{
cat(U)\fr{} \fd{cat(g)} & X \fd{}\ar@/^15pt/[rdd]^{u\circ \phi_1}&\\
cat(V) \ar@/_15pt/[rrd]_{\overline{\phi_2}}\fr{} & \cocartesien Y\ar@{-->}[rd]^{\overline{h}} &\\
&&\overline{Z}}
\]
One obtains the commutative diagram
\[
\xymatrix{ && cat(U) \ar@{->}[dl]_{cat(g)}\fr{}\fd{cat(g)} & X
\ar@{->} `u[l]  `[lll] `[dd]_{\phi_1} [lldd]
\fd{}\\
&cat(V)\!\restriction_{g(X^0)} \fd{\phi_2}\ar@{^{(}->}[r]& cat(V) \fd{\overline{\phi_2}} \fr{} & \cocartesien Y\ar@{->}[ld]^{\overline{h}}\\
&Z \fr{u} & \cocartesien \overline{Z}&}
\]
So $h=\overline{h}\!\restriction_{f(X^0)}$ makes the following
diagram commutative:
\[
\xymatrix{
cat(U)\fr{} \fd{cat(g)} & X \fd{}\ar@/^15pt/[rdd]^{\phi_1}&\\
cat(V)\!\restriction_{g(X^0)} \ar@/_15pt/[rrd]_{\phi_2}\fr{} & Y\!\restriction_{f(X^0)}\ar@{-->}[rd]^{h} &\\
&&Z}
\]
Therefore the following square of flows is cocartesian:
\[
\xymatrix{
cat(U) \fr{}\fd{cat(g)} & X \fd{}\\
cat(V)\!\restriction_{g(X^0)} \fr{} & \cocartesien
Y\!\restriction_{f(X^0)}}
\]
Since $cat(U)\iso cat(V)\!\restriction_{g(X^0)}$ with the first
part of the proof, one gets $X\iso Y\!\restriction_{f(X^0)}$. \epf

\bth \label{singleton} 
Let $f$ be a morphism of flows from $X$ to $Y$. Assume that $f$ is the
pushout of a morphism of flows of the form
$cat(g):cat(U)\longrightarrow cat(V)$ where $g:U\longrightarrow V$ is
a T-homotopy equivalence of globular complexes. For any $\alpha\in Y^0
\backslash f(X^0)$, the topological spaces
${\P}_\alpha^- Y$ and ${\P}_\alpha^+ Y$ are singletons. 
\eth

\bpf Let us suppose first that $X=cat(U)$, $Y=cat(V)$ and $f=cat(g)$
for some T-homotopy equivalence of globular complexes
$g:U\longrightarrow V$. Let $\alpha\in Y^0 \backslash f(X^0)$. One
sees by induction on the globular decomposition of $U$ that the
topological spaces ${\P}^\pm_{g^{-1}(\alpha)}X$ are singletons.  Since
one has ${\P}^\pm_{\alpha}Y\iso {\P}^\pm_{g^{-1}(\alpha)}X$ as
topological spaces, the proof is complete in that case.

Let us take now a general T-homotopy equivalence of flows $h$ from $X$
to $Y$. By hypothesis, there exists a cocartesian diagram of flows
\[
\xymatrix{
cat(U) \fr{}\fd{cat(g)} & X \fd{h}\\
cat(V) \fr{\phi} & \cocartesien Y}
\]
for some T-homotopy equivalence of globular complexes $g:U\longrightarrow V$.
Let $\alpha\in Y^0\backslash h(X^0)$. Since one has the
cocartesian diagram of sets
\[
\xymatrix{
U^0 \fr{}\fd{g^0} & X^0 \fd{h^0}\\
V^0 \fr{\phi} & \cocartesien Y^0}
\]
then there exists a unique $\beta\in V^0\backslash U^0$ such that
$\phi(\beta)=\alpha$. By the first part of this proof, both
topological spaces $ {\P}^\pm_\beta V$ are singletons. Let
$\gamma\in {\P} Y$ with $s(\gamma)=\alpha$. Then one has 
$\gamma=\gamma_1*\dots *\gamma_n$ where the $\gamma_i$ are either
execution paths of ${\P} V$ or execution paths of ${\P} X$. Since
$\alpha=s(\gamma_1)$ and since $\alpha\in Y^0\backslash h(
X^0)$, one deduces that  $\gamma_1\in {\P} V$. But since $\gamma$ is
$\Rm$-equivalent to $\gamma_1$, one deduces that ${\P}^-_\alpha Y$
is a singleton. In the same way, one can check that ${\P}^+_\alpha
Y$ is a singleton as well. \epf

\bth\label{entoure} Let $f$ be a morphism of flows from $X$ to
$Y$. Assume that $f$ is the pushout of a morphism of flows of the form
$cat(g):cat(U)\longrightarrow cat(V)$ where $g:U\longrightarrow V$ is
a T-homotopy equivalence of globular complexes. For any $\alpha\in
Y^0\backslash f(X^0)$, there are execution paths $u$ and $v$ in $Y$
such that $s(u)\in f^0(X^0)$, $t(u)=y$, $s(v)=y$ and $t(v)\in
f^0(X^0)$. \eth

\bpf First suppose that $X=cat(U)$, $Y=cat(V)$ and $f=cat(g)$ for
some T-homotopy equivalence of globular complexes $g:U\longrightarrow
V$. Let $\alpha\in Y^0\backslash g(X^0)$. Then $g^{-1}(\alpha)$ is in
the middle of a globe of the globular decomposition of $X$. In other
terms, there exists $\gamma\in {\P}^{top}X$ such that $\alpha\in
\gamma( ]0,1[)$. So there exists $\gamma_1\in {\P} Y$ and
$\gamma_2\in {\P} Y$ such that $s(\gamma_1)\in g(X^0)$,
$t(\gamma_2)\in g(X^0)$ and $t(\gamma_1)=s(\gamma_2)=\alpha$.
Hence the conclusion in that case.

Take now a general T-homotopy equivalence of flows $h$ from $X$ to
$Y$. By hypothesis, there exists a cocartesian diagram of flows
\[
\xymatrix{
cat(U) \fr{}\fd{cat(g)} & X \fd{h}\\
cat(V) \fr{\phi} & \cocartesien Y}
\]
for some T-homotopy equivalence of globular complexes $g:U\longrightarrow V$.
Let $\alpha\in Y^0\backslash h(X^0)$. Like in
Theorem~\ref{singleton}, there exists a unique $\beta\in
V^0\backslash U^0$ such that $\phi(\beta)=\alpha$. Then using the
first part of this proof, there exist $\gamma_1\in {\P} V$ and
$\gamma_2\in {\P} V$ such that $s(\gamma_1)\in f(U^0)$,
$t(\gamma_2)\in f(U^0)$ and $t(\gamma_1)=s(\gamma_2)=\beta$. Then
$s(\phi(\gamma_1))\in h(X^0)$, $t(\phi(\gamma_2))\in h(X^0)$ and
$t(\phi(\gamma_1))=s(\phi(\gamma_2))=\alpha$. Hence the
conclusion in the general case. \epf

\begin{cor} \label{precorT} Let $f$ be a morphism of flows from $X$ to $Y$. Assume that
$f$ is the pushout of a morphism of flows of the form
$cat(g):cat(U)\longrightarrow cat(V)$ where $g:U\longrightarrow V$ is a T-homotopy
equivalence of globular complexes. Then $f$ is a T-homotopy
equivalence of flows. \end{cor}

\bpf This is an immediate consequence of Theorem~\ref{injecte},
Theorem~\ref{singleton} and Theorem~\ref{entoure}. \epf

\subsection{Comparison with T-homotopy of globular complexes}

\bth\label{corT} Let $X$ and $U$ be globular complexes. 
Let $f:X\longrightarrow U$ be a T-homotopy equivalence of globular
complexes. Then $cat(f):cat(X)\longrightarrow cat(U)$ is a T-homotopy
equivalence of flows. Conversely, if $g:cat(X)\longrightarrow cat(U)$
is a T-homotopy equivalence of flows, then $g=cat(f)$ for some
T-homotopy equivalence $f:X\longrightarrow U$ of globular complexes.
\eth

\bpf Let $f:X\longrightarrow U$ be a T-homotopy equivalence of globular
complexes. Then $cat(f):cat(X)\longrightarrow cat(U)$ is a T-homotopy
equivalence of flows by Corollary~\ref{precorT}.

Conversely, let $X$ and $U$ be two globular complexes. Let
$g:cat(X)\longrightarrow cat(U)$ be a T-homotopy equivalence of
flows. Let $(\de Z_\beta, Z_\beta,\phi_\beta)_{\beta<\lambda}$ be the
globular decomposition of $X$. The morphism $g$ gives rise to a
one-to-one set map $g^0$ from $cat(X)^0$ to $cat(U)^0$ and to an
homeomorphism ${\P}g:{\P} cat(X)\longrightarrow{\P}
cat(U)\!\restriction_{g^0(cat(X)^0)}$. Let $i_U:{\P} U\to {\P}^{top}U$
given by Theorem~\ref{i}.

Let us suppose by induction on $\beta$ that there exists a one-to-one
morphism of globular complexes $f_\beta:X_\beta\longrightarrow U$ such
that ${\P} f_\beta:{\P} X_\beta\longrightarrow {\P} U$ coincides with
the restriction of ${\P}g$ to ${\P} X_\beta$. One has to prove that
the same thing holds for $\beta+1$. There is a cocartesian diagram of
multipointed topological spaces
\[
\xymatrix{
\glob^{top}(\de Z_\beta)\fd{}\fr{\phi_\beta} & X_\beta \fd{}\\
\glob^{top}(Z_\beta)\fr{\overline{\phi_\beta}} & \cocartesien
X_{\beta+1}}
\]
Let \[k(z,t)= i_U(g\circ
\overline{\phi_\beta}(z))(t)\] for
$z\in Z_\beta$ and $t\in [0,1]$. For any $z\in Z_\beta$, $k(z,-)$
is an execution path of $U$. The composite of morphisms of
globular complexes
\[\xymatrix{\glob^{top}(\de Z_\beta)\fr{\phi_\beta}&X_\beta\fr{f_\beta}& U}\]
gives rise to an execution path $\ell(z,t)$ for any $z\in \de
Z_\beta$.  Since $\de Z_\beta$ is compact, then 
there exists a continuous map $\psi:\de
Z_\beta\p [0,1]\longrightarrow [0,1]$ such that
$\ell(z,t)=k(z,\psi(z,t))$ for any $z\in \de Z_\beta$ and any
$t\in [0,1]$ by Proposition~\ref{span}. Therefore the mapping
\[\overline{k}:(z,t)\mapsto k(z,\mu(z)t+(1-\mu(z))\psi(z,t))\]
induces a morphism of globular complexes
$f_{\beta+1}:X_{\beta+1}\longrightarrow U$ which is an extension
of $f_{\beta}:X_{\beta}\longrightarrow U$.

One now wants to prove that the restriction of $f_{\beta+1}$ to
$\glob^{top}(Z_\beta\backslash \de Z_\beta)$ is one-to-one.  Suppose
that there exists two points $(z,t)$ and $(z',t')$ of
$\glob^{top}(Z_\beta\backslash \de Z_\beta)$ such that
$f_{\beta+1}(z,t)=f_{\beta+1}(z',t')$. If $z\neq z'$, then $g\circ
\overline{\phi_\beta}(z)\neq g\circ
\overline{\phi_\beta}(z')$ since $g$ is one-to-one. So the two
execution paths $\overline{k}(z,-)$ and $\overline{k}(z',-)$ are two
distinct execution paths intersecting at
$\overline{k}(z,t)=\overline{k}(z',t')$. The latter point necessarily
belongs to $U^0$. Since the topological spaces ${\P}^-_\alpha U$ and
${\P}^+_\alpha U$ are both singletons for $\alpha\in U^0\backslash
\phi(X^0)$, then
\[\overline{k}(z,t)=\overline{k}(z',t')\in \phi(X^0).\]
There are two possibilities:
$\overline{k}(z,t)=\overline{k}(z',t')=g\circ \phi_\beta(0)$ and
$\overline{k}(z,t)=\overline{k}(z',t')=g\circ \phi_\beta(1)$
(notice that $\phi_\beta(0)$ and $\phi_\beta(1)$ may be equal).
The equality $\overline{k}(z,t)=\overline{k}(z',t')=g\circ
\phi_\beta(0)$ implies $t=t'=0$ and the equality
$\overline{k}(z,t)=\overline{k}(z',t')=g\circ \phi_\beta(1)$ implies
$t=t'=1$. In both cases, one has $(z,t)=(z',t')$: contradiction. So
$f_{\beta+1}$ is one-to-one.

If $\beta<\lambda$ is a limit ordinal, then let
$f_\beta=\liminj_{\alpha<\beta}f_\alpha$. The latter map is still a
one-to-one continuous map and a morphism of globular complexes.  So
one obtains a one-to-one morphism of globular complexes
$f:X\longrightarrow U$ such that ${\P} f: {\P} X\longrightarrow {\P}
U$ coincides with ${\P} g$.

Now let us prove that $f$ is surjective. Let $x\in U$. First case:
$x\in U^0$. If $x\notin g(X^0)=f(X^0)$, then by hypothesis,
${\P}_{x}^- cat(U)$ and ${\P}_{x}^+ cat(U)$ are singletons. So $x$
necessarily belongs to an execution path between two points of
$g(X^0)$. Since $g$ is a bijection from ${\P} X$ to ${\P}
U\!\restriction_{g(X^0)}$, this execution path necessarily belongs to
$g(X)$.  Therefore $x\in g(X)$. Second case: $x\in U\backslash
U^0$. Then there exists an execution path $\gamma$ of $U$ passing by
$x$. If $\gamma(0)\notin g(X^0)$ (resp.  $\gamma(1)\notin g(X^0)$),
then there exists an execution path going from a point of $g(X^0)$ to
$\gamma(0)$ (resp. going from $\gamma(1)$ to a point of $g(X^0)$)
because $\gamma(0)$ is not an initial state (resp. a final state) of
$g(X)$.  Therefore one can suppose that $\gamma(0)$ and $\gamma(1)$
belong to $g(X^0)$. Once again we recall that $g$ is a bijection from
${\P} X$ to ${\P} U\!\restriction_{g(X^0)}$, so $x\in f(X)$.
Therefore $U\subset f(X)$. So $f$ is bijective.

At last, one has to check that $f^{-1}:U\longrightarrow X$ is
continuous. Let $T$ be a compact of the globular decomposition of $U$
(not of $X$ !), let $q$ be the corresponding attaching map, and
consider the composite
\[\xymatrix{\glob^{top}(T)\fr{q} & U \fr{f^{-1}}& X}\]
There are then two possibilities.

First of all, assume that $T=\{x\}$ for some $x\in U^0$.  Then there
exists $\gamma,\gamma'\in U$ such that $\gamma*({\P} q)(x)*\gamma'\in
g({\P} X)$. Let $\gamma*{\P} q(x)*\gamma'=g(\gamma'')$. Let $i_X:{\P}
X\longrightarrow {\P}^{top}X$ given by Theorem~\ref{i}. Then the
execution path $i_X(\gamma'')$ of $X$ becomes an execution path
$f\circ i_X(\gamma'')$ of $U$ since $f$ is one-to-one. Let us consider
the execution path $i_U(\gamma*{\P} q(x)*\gamma')$ of $U$. By
Proposition~\ref{span}, there exists a continuous non-decreasing map
$\omega:[0,1]\longrightarrow [0,1]$ such that $\omega(0)=0$,
$\omega(1)=1$ and such that
\[i_U(\gamma*{\P} q(x)*\gamma')=f\circ
i_X(\gamma'')\circ \omega.\] Then $\omega$ is
necessarily bijective, and so an homemorphism since $[0,1]$ is
compact.  Therefore  $f^{-1}\circ i_U(\gamma*{\P}
q(x)*\gamma')=i_X(\gamma'')\circ
\omega $.  So  $f^{-1}\circ
q(\glob^{top}(T))$ is a compact of $X$.

Now suppose that $T$ contains more than one element. Then
${\P}_{q(0)}^-U$ and ${\P}_{q(1)}^+U$ are not singletons.  So
$q(0)$ and $q(1)$ belong to $g(X^0)=f(X^0)$. Then
${\P}(g^{-1}\circ
q)(T)=(g^{-1}\circ{\P}
q)(T)$ is a compact of ${\P} X$ (since $g$ is an
homeomorphism !). By Proposition~\ref{span}, there exists a
continuous map $\omega:T\longrightarrow \ttop([0,1],[0,1])$ such
that $\omega(0)=0$, $\omega(1)=1$ and such
that $\omega(z)$ is non-decreasing for any $z\in T$
and such that
\[i_U({\P} q(z))=f\circ i_X(g^{-1}\circ {\P} q(z))\circ
\omega(z)\] for any $z\in T$. The map $z\mapsto
i_X(g^{-1}\circ {\P} q(z))\circ
\omega(z)$ is mapped by the set map
\[\top(T,{\P}^{top}X)\longrightarrow \top(T\p [0,1],X)\] to a
function $\omega'\in \top(T\p [0,1],X)$. Therefore
\[f^{-1}\circ q(\glob^{top}(T))=\omega'(T\p
[0,1])\] is again a compact of $X$.

To conclude, let $F$ be a closed subset of $X$. Then
\[f(F)\cap q(\glob^{top}(T)) =
f(F\cap (f^{-1}\circ q)(\glob^{top}(T))).\] Since $f^{-1}\circ
q(\glob^{top}(T))$ is always compact, the set $f(F)\cap
q(\glob^{top}(T))$ is compact as well. Since $U$ is equipped with
the weak topology induced by its globular decomposition, the set
$f(F)$ is a closed subspace of $U$. So $f^{-1}$ is continuous.
\epf

\section{Conclusion}

We have defined in this part a class of morphisms of flows, the
T-homotopy equivalences, such that there exists a T-homotopy
equivalence between two globular complexes if and only if there exists
a T-homotopy equivalence between the corresponding flows.  So not only
the category of flows allows the study of S-homotopy of globular
complexes, but also the study of T-homotopy of globular complexes.

\part{Application : the underlying homotopy type of a flow}\label{uhtf}

\section{Introduction}

The main theorem of this paper (Theorem~\ref{equiS}) establishes the
equivalence of two approaches of dihomotopy. The first one uses the
category of globular complexes in which the concurrent processes are
modelled by topological spaces equipped with an additional structure,
the globular decomposition, encoding the time flow, and in which the
execution paths are "locally strictly increasing" continuous maps.
The second one uses the category of flows in which the concurrent
processes are modelled by categorical-like objects and in which it is
possible to define a model structure relevant for the study of
dihomotopy.

Another interest of this equivalence is that it makes the construction
of the \textit{underlying homotopy type} of a flow possible.  Indeed,
loosely speaking, a dihomotopy type is an homotopy type equipped with
an additional structure encoding the time flow. So there must exist a
forgetful functor $|-|:\dtop\longrightarrow \ho(\top)$ from the
category of flows to the category of homotopy types which is also a
dihomotopy invariant, i.e. sending weak S-homotopy and T-homotopy
equivalences to isomorphisms.

\section{Construction of the underlying homotopy type functor}

\bd (cf. Part~\ref{part1} Section~\ref{def}) 
Let $(X,X^0)$ be a multipointed topological space. Then the mapping
\[(X,X^0)\mapsto X\] induces a functor $|-|:\top^m\longrightarrow \top$
called the {\rm underlying topological space} of $(X,X^0)$. \ed

\bp The underlying topological space construction induces a functor
$|-|:\glCW\longrightarrow \top$ from the category of globular
CW-complexes to the category of topological spaces. Moreover, for any
S-homotopy equivalence $f:X\longrightarrow U$ of globular
CW-complexes, the continuous map $|f|:|X|\longrightarrow |U|$ is a
homotopy equivalence of topological spaces. \ep

\bpf It suffices to prove that if $f$ and $g$ are two morphisms
of globular complexes which are S-homotopy equivalent, then $|f|$ and
$|g|$ are two homotopy equivalent continuous maps. Let $H$ be a
S-homotopy between $f$ and $g$. By Proposition~\ref{caracgl}, $H$
induces a continuous map $h\in\top([0,1],\tgltop(X,Y))$, so a
continuous map $h\in\top([0,1],\ttop(|X|,|Y|))$. Hence an homotopy
between the continuous maps $|f|$ and $|g|$. \epf

\begin{cor} The functor $|-|:\glCW\longrightarrow \top$ induces a unique
functor $|-|:\glCW[\mathcal{SH}^{-1}]\longrightarrow \ho(\top)$
making the following diagram commutative:
\[
\xymatrix{
\glCW \fd{} \fr{} & \ho(\top) \\
\glCW[\mathcal{SH}^{-1}]\ar@{->}[ru] & }
\]
\end{cor}

\bd The composite functor
\[\xymatrix{
|-|:\dtop\fr{}& \dtop[\mathcal{S}^{-1}] \simeq
\glCW[\mathcal{SH}^{-1}]\fr{|-|}& \ho(\top)}
\]
is called the {\rm underlying homotopy type functor}. If $X$ is a
flow, then $|X|$ is called the {\rm underlying homotopy type} of
$X$. \ed

\bp If $f:X\longrightarrow Y$ is a weak S-homotopy equivalence of
flows, then $|f|$ is an isomorphism of $\ho(\top)$. If
$g:cat(X)\longrightarrow cat(Y)$ is a T-homotopy equivalence of
flows, then $|g|$ is an isomorphism of $\ho(\top)$ as well. \ep

\bpf Obvious with Theorem~\ref{corT}. \epf

Figure~\ref{ex1} represents the simplest example of T-homotopy
equivalence. The underlying homotopy types of its source and its
target are both equal to the homotopy type of the point.

Notice that the functor from $\dtop$ to $\ho(\top)$ defined by
associating to a flow $X$ the homotopy type of the disjoint sum $\P
X\sqcup X^0$ is not a dihomotopy invariant. Therefore the functor
$X\mapsto \hbox{"homotopy type of }\P X\sqcup X^0\hbox{"}$ has no
relation with the underlying homotopy type functor. In the case of
Figure~\ref{ex1}, the discrete space $\{u,s(u),t(u)\}$ becomes the
discrete space $\{v,w,v*w,s(v),t(w),t(v)=s(w)\}$.

\begin{question} How to define the underlying homotopy type of a
flow without using the category of globular complexes ?
\end{question}

\section{Conclusion}

The underlying homotopy type functor is a new dihomotopy invariant
which can be useful for the study of flows up to dihomotopy.

\ifx\undefined\bysame
\newcommand{\bysame}{\leavevmode\hbox to3em{\hrulefill}\,}
\fi

\end{document}